\documentstyle[10pt,amscd,xypic,amssymb,combelow,enumitem]{amsart}

\xyoption{all}
\CompileMatrices
\newcommand{\nc}{\newcommand}

\makeatletter
\@addtoreset{equation}{section}
\makeatother

\newenvironment{proof}{{\noindent \textbf{Proof}\,\,}}{\hspace*{\fill}$\Box$\medskip}

\newtheorem{theorem}[equation]{Theorem}
\newtheorem{proposition}[equation]{Proposition}

\theoremstyle{definition}

\newtheorem{definition}[equation]{Definition}

\theoremstyle{remark}

\newtheorem{remark}[equation]{Remark}

\nc{\fa}{{\mathfrak{a}}}
\nc{\fb}{{\mathfrak{b}}}
\nc{\fg}{{\mathfrak{g}}}
\nc{\fh}{{\mathfrak{h}}}
\nc{\fj}{{\mathfrak{j}}}
\nc{\fn}{{\mathfrak{n}}}
\nc{\fu}{{\mathfrak{u}}}
\nc{\fp}{{\mathfrak{p}}}
\nc{\fr}{{\mathfrak{r}}}
\nc{\ft}{{\mathfrak{t}}}
\nc{\fsl}{{\mathfrak{sl}}}
\nc{\fgl}{{\mathfrak{gl}}}
\nc{\hsl}{{\widehat{\mathfrak{sl}}}}
\nc{\hgl}{{\widehat{\mathfrak{gl}}}}

\nc{\fC}{{\mathfrak{C}}}
\nc{\fZ}{{\mathfrak{Z}}}
\nc{\fW}{{\mathfrak{W}}}

\nc{\pol}{{\text{Poles}}}

\nc{\BA}{{\mathbb{A}}}
\nc{\BC}{{\mathbb{C}}}
\nc{\BK}{{\mathbb{K}}}
\nc{\BM}{{\mathbb{M}}}
\nc{\BN}{{\mathbb{N}}}
\nc{\BF}{{\mathbb{F}}}
\nc{\BQ}{{\mathbb{Q}}}
\nc{\BP}{{\mathbb{P}}}
\nc{\BR}{{\mathbb{R}}}
\nc{\BZ}{{\mathbb{Z}}}

\nc{\CA}{{\mathcal{A}}}
\nc{\CB}{{\mathcal{B}}}
\nc{\DD}{{\mathcal{D}}}
\nc{\CE}{{\mathcal{E}}}
\nc{\CF}{{\mathcal{F}}}
\nc{\tCF}{{\widetilde{\CF}}}
\nc{\oCF}{{\overline{\CF}}}
\nc{\CG}{{\mathcal{G}}}
\nc{\CL}{{\mathcal{L}}}
\nc{\CM}{{\mathcal{M}}}
\nc{\CH}{{\mathcal{H}}}
\nc{\CN}{{\mathcal{N}}}
\nc{\CO}{{\mathcal{O}}}
\nc{\CP}{{\mathcal{P}}}
\nc{\CQ}{{\mathcal{Q}}}
\nc{\CR}{{\mathcal{R}}}
\nc{\CS}{{\mathcal{S}}}
\nc{\CT}{{\mathcal{T}}}
\nc{\CU}{{\mathcal{U}}}
\nc{\CV}{{\mathcal{V}}}
\nc{\CW}{{\mathcal{W}}}
\nc{\tCW}{{\widetilde{\CW}}}
\nc{\oCW}{{\overline{\CW}}}

\nc{\uu}{{U_q(\hgl_n)}}
\nc{\uuo}{{U^0_q(\hgl_n)}}
\nc{\uug}{{U_q^+(\hgl_n)}}
\nc{\uul}{{U_q^-(\hgl_n)}}

\nc\bara{{\bar{a}}}
\nc\bari{{\bar{i}}}
\nc\barj{{\bar{j}}}
\nc\bark{{\bar{k}}}
\nc\barii{{\bar{i'}}}
\nc\barjj{{\bar{j'}}}

\nc{\su}{{U_q(\hsl_n)}}
\nc{\suo}{{U^0_q(\hsl_n)}}
\nc{\sug}{{U_q^+(\hsl_n)}}
\nc{\sul}{{U_q^-(\hsl_n)}}

\nc{\uui}{{U_q(\hgl_1)}}
\nc{\uuio}{{U^0_q(\hgl_1)}}
\nc{\uuig}{{U_q^+(\hgl_1)}}
\nc{\uuil}{{U_q^-(\hgl_1)}}

\nc{\e}{{\varepsilon}}
\nc{\nn}{{\mathbb{N}}^n}
\nc{\Tr}{{\text{Tr}}}
\nc{\tT}{{T}}

\nc{\od}{{\overline{d}}}
\nc{\rg}{{\textrm{R}\Gamma}}
\nc{\erg}{{\emph{R}\Gamma}}
\nc{\id}{{\textrm{id}}}

\def\bd{{\mathbf{d}}}

\def\bk{{\mathbf{k}}}
\def\bl{{\mathbf{l}}}

\def\bde{{\boldsymbol{\delta}}}
\def\bs{{\boldsymbol{\varsigma}}}

\def\pt{\textrm{pt}}

\def\ofZ{{\overline{\fZ}}}

\def\ev{\textrm{ev}}

\def\ev{\textrm{ev}}

\def\orr{{\tilde{r}}}
\def\oe{{\tilde{e}}}
\def\of{{\tilde{f}}}

\def\op{{\overline{p}}}
\def\oq{{\overline{q}}}

\def\tA{{\tilde{A}}}

\def\tZ{{\tilde{Z}}}
\def\bx{{\textbf{x}}}

\def\id{{\text{Id}}}
\def\col{{\text{col }}}
\def\sym{{\text{Sym}}}

\def\Tr{{\text{Tr}}}

\nc{\loccitt}{{\emph{loc. cit.}}}
\nc{\loccit}{{\emph{loc. cit. }}}

\begin{document}

\title[Affine Laumon Spaces and Integrable Systems]{\Large{\textbf{Affine Laumon Spaces and Integrable Systems}}}

\author[Andrei Negu\cb t]{Andrei Negu\cb t}
\address{MIT, Department of Mathematics, Cambridge, MA, USA}
\address{Simion Stoilow Institute of Mathematics, Bucharest, Romania}
\email{andrei.negut@@gmail.com}

\maketitle

\begin{abstract}

In this paper we formalize and prove a conjecture of Braverman (\cite{B}) concerning integrals of the Chern polynomial of the tangent bundle to affine Laumon spaces. This provides the computation of the Nekrasov partition function (\cite{N}) of ${\mathcal{N}} = 2$ gauge theory with adjoint matter on $\BC^2$ in the $\Omega$--background, in the presence of a full surface operator insertion. 

\end{abstract}

\section{Introduction}

\noindent Affine Laumon spaces\footnote{The name ``affine Laumon spaces" is short-hand for the phrase ``certain quasiprojective varieties whose representation theory is controlled by $\hgl_n$ in the same way that the representation theory of Laumon's compactification of $\text{Maps}(\BP^1, \text{flag variety})$ is controlled by $\fgl_n$"} appear naturally in geometic representation theory as semismall resolutions of singularities of Uhlenbeck spaces for the affine Lie algebra $\hgl_n$. Uhlenbeck spaces appear in gauge theory as partial compactifications of moduli spaces of $U(n)-$instantons, hence the partition function of pure gauge theory on $\BC^2$ (in the $\Omega$--background and with a full surface operator insertion) is given by:
$$
\int_{\text{Uhlenbeck space}} 1 = \int_{\text{affine Laumon space}} 1
$$
(see \cite{B2} for a mathematical treatment of the notions above) where the latter equality follows because affine Laumon spaces are resolutions of singularities of Uhlenbeck spaces. To introduce adjoint matter into the gauge theory, the singular Uhlenbeck spaces are not the correct algebro-geometric object. Instead, one needs to replace 1 by the Chern polynomial of the tangent bundle to the smooth affine Laumon spaces, hence the partition function of gauge theory with adjoint matter (in the $\Omega$--background and with a full surface operator insertion, see \cite{N}) is:
\begin{equation}
\label{eqn:def z}
\tZ_m = \sum_{\bd \in \nn} \bx^{\bd} \int_{\CM_{\bd}} c (T\CM_{\bd}, 2m)
\end{equation}
where $\{\CM_{\bd}\}_{\bd \in \nn}$ denote the connected components of affine Laumon space, which are indexed by $\bd = (d_1,...,d_n) \in \nn$, and $\bx^{\bd} = x_1^{d_1}...x_n^{d_n}$ are formal variables. The parameter $m$ tracks the mass of the adjoint matter. All integrals in this paper are equivariant with respect to the torus $T = (\BC^*)^n \times \BC^* \times \BC^*$ which acts on $\CM_{\bd}$, whose equivariant parameters will be denoted by $a_1,...,a_n, 1, p$. Therefore \eqref{eqn:def z} is a power series in $x_1,...,x_n$ with coefficients which are rational functions in $a_1,...,a_n,p$. \\

\noindent In \cite{B}, Braverman asked to compute the generating function \eqref{eqn:def z} by relating it to a non-stationary deformation of the affine trigonometric Calogero-Moser hamiltonian: 
\begin{equation}
\label{eqn:cal}
\CH_m = \sum_{i=1}^n D_i \left(D_i - D_{i-1} + \frac pn \right) + \sum_{(i,j) \in \frac {\BZ^2}{(n,n)\BZ}}^{i<j} \frac {m(m+1)}{\left(1 - x_i ... x_{j-1} \right)\left(1 - x_i^{-1} ... x_{j-1}^{-1} \right)} \qquad 
\end{equation}
where $p+n$ is the central charge of $\hgl_n$, and we set $D_i = x_i \partial_{x_i}$ and $D_0 = D_n$.\footnote{We refer to (6.1) of \cite{EK} for the definition of the hamiltonian \eqref{eqn:cal}, and references therein for the history of this differential operator. The differences between our operator and that of \loccit is that our corresponds to $\hgl_n$, whereas theirs corresponds to $\hsl_n$} For generic parameters $b_1,...,b_n$, the differential operator $\CH$ has a unique eigenfunction:
\begin{equation}
\label{eqn:eigen}
Y_m \in x_1^{b_1}...x_n^{b_n} \Big(1 + \BC(b_1,...,b_n,m,p)[[x_1,...,x_n]] \Big) 
\end{equation}
namely:
\begin{equation}
\label{eqn:function}
\CH_m \cdot Y_m = \sum_{i=1}^n b_i \left(b_i - b_{i-1} + \frac pn \right)Y_m
\end{equation}
In \cite{Laumon}, we showed that the functions $\tZ_m$ and $Y_m$ are closely related in the special case of the usual Laumon spaces, which corresponds to setting $x_n = 0$ in all formulas above. In the present paper, we prove the general relationship: \\

\begin{theorem}
\label{thm:main}

The generating function \eqref{eqn:def z} is given by:
\begin{equation}
\label{eqn:finally equal}
\tZ_m = \frac{Y_m}{x_1^{b_1}...x_n^{b_n} \cdot \delta^{m+1}}
\end{equation}
where the formal series:
\begin{equation}
\label{eqn:weyl}
\delta = \prod_{(i,j) \in \frac {\BZ^2}{(n,n)\BZ}}^{i<j} (1-x_i...x_{j-1})
\end{equation}
is called the Weyl determinant of $\hgl_n$, and the equivariant parameters $a_1,...,a_n$ that define $\tZ_m$ are connected to the parameters $b_1,...,b_n = b_0$ that define $Y_m$ by:
\begin{equation}
\label{eqn:connection}
b_i - b_{i-1} = a_i - \frac {a_1+...+a_n}n + \frac pn \left(i - \frac {n+1}2 \right) \qquad \forall i \in \{1,...,n\}
\end{equation}

\end{theorem}

\noindent As in \cite{Laumon}, the main idea of the proof is to present the function $\tZ_m$ as the character of the so-called Ext operator, inspired by the construction of \cite{CO} for Hilbert schemes:
\begin{equation}
\label{eqn:op cohom}
\tA_m : H \rightarrow H \qquad \text{ where } \qquad H = \bigoplus_{\bd \in \nn} H^*_T(\CM_{\bd})
\end{equation}
The main geometric computation performed in the present paper relates the operator \eqref{eqn:op cohom} to intertwiners of the action $\hgl_n \curvearrowright H$ that was defined in \cite{Aff} (generalizing the action $\hsl_n \curvearrowright H$ of \cite{FK}). Then we can use the tools of \cite{EK} to relate the character of $\hgl_n$ intertwiners with the eigenfunction \eqref{eqn:eigen}, thus proving Theorem \ref{thm:main}. \\

\noindent Equivariant cohomology has a natural deformation, namely equivariant $K$--theory. In this language, the generating function \eqref{eqn:def z} admits the following deformation:
\begin{equation}
\label{eqn:def z 2}
Z_m = \sum_{\bd \in \nn} \bx^{\bd} \cdot \chi_T \left( {\CM_{\bd}},  \sum_{k=0}^\infty (-q^{2m})^k \cdot \wedge^k(T\CM^\vee_{\bd}) \right)
\end{equation}
Note that the coefficients of the generating series above are rational functions in the equivariant parameters $u_1,...,u_n,q,\oq$ of $T$. By letting $u_i = e^{\hbar a_i}$, $q = e^{\hbar}$, $\oq = e^{\hbar p}$ and taking the leading order term as $\hbar \rightarrow 0$, the expression \eqref{eqn:def z 2}  degenerates to \eqref{eqn:def z}. From a physical point of view, the generating function \eqref{eqn:def z 2} is the partition function of ${\mathcal{N}} = 1$ supersymmetric gauge theory on $\BC^2$ (in the $\Omega$--background and with a full surface operator insertion) times a small circle. \\

\noindent One can present $Z_m$ as the character of the $K$--theoretic version of the Ext operator:
\begin{equation}
\label{eqn:op k}
A_m : K \rightarrow K \qquad \text{ where } \qquad K = \bigoplus_{\bd \in \nn} K_T(\CM_{\bd})
\end{equation}
The algebra $\uu$ acts on $K$ (see \cite{Aff}, generalizing the action of $\su$ from \cite{BF}, \cite{FFNR}, \cite{T}). Moreover, we show that the operator $A_m$ is uniquely determined by the way it interacts with the generators of the quantum group, namely Propositions \ref{prop:comm root e}, \ref{prop:comm root f}, \ref{prop:comm group}. These properties determine the operator \eqref{eqn:op k} completely, and can be interpreted as a mathematical incarnation of the AGT correspondence with bifundamental matter and a full surface operator insertion (in the case at hand, the conformal field theory side of AGT is provided by the quantum group $\uu$, which is nothing but the $_qW$--algebra of $\fgl_n$ corresponding to the zero nilpotent). \\

\noindent However, if one were able to relate the operator $A_m$ to actual $\uu$ intertwiners, then one could use the results of \cite{ESV} and \cite{S} to deduce that the latter intertwiners are eigenfunctions of certain difference operators. This would lead to an analogue of Theorem \ref{thm:main} for $K$--theory instead of cohomology, or in other words, it would yield a computation of \eqref{eqn:def z 2} instead of \eqref{eqn:def z}. On a different note, when $x_n = 0$ the function \eqref{eqn:def z 2} was computed in \cite{BFS}, but their proof does not study the operator $A_m$. \\

\noindent I would like to thank Andrei Okounkov for first introducing me to this problem, back in 2007, as well as for all his guidance and help along the way. I would like to thank Alexander Braverman for initially posing this question, as well as for all the interest and explanations he provided along the years. I would also like to thank Pavel Etingof, Michael Finkelberg, Dennis Gaitsgory, Davesh Maulik, Valerio Toledano-Laredo and Alexander Tsymbaliuk for a great host of useful discussions, as well as Sachin Gautam for his hospitality while this paper was being written. I gratefully acknowledge the support of NSF grant DMS--1600375. \\

\section{Quantum groups and the $K$--theory of affine Laumon spaces} 
\label{sec:affine}

\subsection{}\label{sub:intro}

Let $n \geq 2$. We will work with the monoid $\nn$, whose elements will be called {\bf degree vectors} and will be denoted by $\bk = (k_1,...,k_n)$. Consider the pairing:
$$
\langle \cdot , \cdot \rangle : \nn \times \nn \longrightarrow \BZ \qquad \qquad \langle \bk, \bl \rangle = \sum_{i=1}^{n} k_i l_i - k_i l_{i+1}
$$
where we identify $l_{n+1} = l_1$. We also consider the standard generators of $\nn$:
$$
\bs^i = (\underbrace{0,...,0,1,0,...,0}_{i-\text{th position}})
$$
and observe that:
$$
\bde = \bs^1 +... +\bs^{n} = (1,...,1)
$$
spans the kernel of the bilinear form. We will write $\bk+1$ for the degree vector whose $i-$th position is $k_{i-1}$. This convention is chosen such that $\bs^{i}+1 = \bs^{i+1}$. Finally, for any additive variables $a_1,...,a_n$, we set:
\begin{equation}
\label{eqn:additive}
a_\bk = \sum_{i=1}^{n} a_i k_i \qquad \forall \bk = (k_1,...,k_n) \in \nn
\end{equation}
and similarly for multiplicative variables $u_1,...,u_n$, we set $u_\bk = \prod_{i=1}^{n} u_i^{k_i}$. \\


\subsection{}\label{sub:quantum1}


The {\bf quantum group} associated to the $n$ vertex cyclic quiver is the algebra:
$$
\su  := \BC(q) \Big \langle e_1,...,e_n,f_1,...,f_n, q^{\pm h_1}...,q^{\pm h_n}, q^\gamma \Big \rangle
$$ 
modulo the following relations, for all indices $i,j$ modulo $n$:
\begin{align}
&q^{h_j} e_i q^{-h_j} = q^{+ \langle \bs^i, \bs^j \rangle} e_i \label{eqn:sug1} \\
&q^{h_j} f_i q^{-h_j} = q^{- \langle \bs^i, \bs^j \rangle} f_i \label{eqn:sug1 prim}
\end{align}
\begin{equation}
\label{eqn:sug2}
\left [e_i, f_j \right ] = \delta_j^i \cdot \frac {q^{h_i-h_{i+1}} - q^{h_{i+1}-h_i}}{q-q^{-1}} 
\end{equation}
as well as:
\begin{equation}
\label{eqn:sug3}
[e_i, e_j] = 0, \ \text{if } i - j \not \equiv \{-1,1\} \text{ mod } n
\end{equation}
\begin{equation}
\label{eqn:sug4}
\left[ e_i, [e_i, e_{i + 1}]_q \right]_{\frac 1q} = \left[ e_i, [e_i, e_{i-1}]_q \right]_{\frac 1q} = 0
\end{equation}
where $[x,x']_q = xx'-qx'x$, together with the analogous relations for $f$'s instead of $e$'s. The element $\gamma$ is central, and we will write $e_{i+n} = e_i$, $f_{i+n} = f_i$ and $h_{i+n} = h_i - \gamma$ for all integers $i$. The assignments: 
\begin{equation}
\label{eqn:cop1}
\Delta(h_i) = h_i \otimes 1 + 1 \otimes h_i
\end{equation}
\begin{equation}
\label{eqn:cop2}
\Delta(e_i) = q^{h_i-h_{i+1}} \otimes e_i + e_i \otimes 1
\end{equation}
\begin{equation}
\label{eqn:cop3}
\Delta(f_i)  =  1 \otimes f_i + f_i \otimes q^{h_{i+1} - h_i}
\end{equation}
for all $i\in \BZ$ give rise to a bialgebra structure on $\su$. \\

\subsection{}\label{sub:quantum2}

Consider the {\bf quantum Heisenberg algebra}:
$$
\uui := \BC(q) \Big \langle p_{\pm 1}, p_{\pm 2}, p_{\pm 3}, ... , q^{\pm \gamma} \Big \rangle
$$
where $\gamma$ is a central element and the generators $p_{k}$ satisfy the relations:
\begin{equation}
\label{eqn:heis rel}
[p_k, p_l] = \delta_{k+l}^0 k \cdot \text{constant}_k \cdot (q^{k\gamma} - q^{-k\gamma})  
\end{equation}
where ``$\text{constant}_k$" denotes an arbitrary element of the ground field, although for the purposes of the present paper, we will fix a certain choice in \eqref{eqn:fix constant}. The assignments $\Delta(\gamma) = \gamma \otimes 1 + 1 \otimes \gamma$ and:
\begin{equation}
\label{eqn:prim1}
\Delta(p_k) = q^{k\gamma} \otimes p_k + p_k \otimes 1 
\end{equation}
\begin{equation}
\label{eqn:prim2}
\Delta(p_{-k}) = 1 \otimes p_{-k} + p_{-k} \otimes q^{-\gamma k}
\end{equation}
for all $k\in \BN$ give rise to a bialgebra structure on $\uui$. In the present paper, we will work with a related set of generators of the quantum Heisenberg algebra:
\begin{equation}
\label{eqn:group-like}
\sum_{k=0}^\infty g_{\pm k} z^k := \exp \left(\sum_{k=1}^\infty \frac {p_{\pm k} z^k}k \right) 
\end{equation}
It is easy to see that \eqref{eqn:prim1} and \eqref{eqn:prim2} imply the following coproduct relations:
\begin{equation}
\label{eqn:group1}
\Delta(g_k) = \sum_{a+b = k} g_a q^{b\gamma} \otimes g_b
\end{equation}
\begin{equation}
\label{eqn:group2}
\Delta(g_{-k}) = \sum_{a+b = k} g_a \otimes g_b q^{-a\gamma}
\end{equation}

\subsection{}\label{sub:quantum3} 

The {\bf affine quantum group} that will feature in the present paper is:
\begin{equation}
\label{eqn:quantum}
\uu = \su \otimes \uui
\end{equation}
In \cite{Tor}, we worked with the equivalent RTT presentation of the quantum group (see \cite{FRT}, \cite{RS}, \cite{DF}, \cite{GM}), which entails the existence of a family of {\bf root generators}:
\begin{equation}
\label{eqn:root}
e_{[i;j)}, \ f_{[i;j)} \in \uu
\end{equation}
$\forall i < j$, such that $e_i = e_{[i;i+1)}$ and $f_i = f_{[i;i+1)}$. The root generators satisfy the following relations (we write $\bari$ for the residue of the integer $i$ in the set $\{1,...,n\}$):
\begin{align}
&q^{h_k} e_{[i;j)} q^{-h_k} = q^{\delta_{\bark}^{\bari} - \delta_{\bark}^{\barj}} e_{[i;j)} \label{eqn:rtt cartan 1} \\
&q^{h_k} f_{[i;j)} q^{-h_k} = q^{\delta_{\bark}^{\barj} - \delta_{\bark}^{\bari}} f_{[i;j)} \label{eqn:rtt cartan 2}
\end{align}
and:
$$
\frac {e_{[i;j)} e_{[i';j')}}{q^{\delta_{\bari}^{\barj} - \delta_{\bari}^{\barjj} - \delta_{\barj}^{\barjj}}} - \frac {e_{[i';j')} e_{[i;j)}}{q^{\delta_{\bari}^{\barj} - \delta_{\barj}^{\barii} - \delta_{\barj}^{\barii}}} = \delta_{\barii}^{\barj} e_{[i;j+j'-i')} - \delta_{\bari}^{\barjj} e_{[i+i'-j';j)} + 
$$
\begin{equation}
\label{eqn:rtt 1}
+ (q-q^{-1}) \left(\sum_{i' < a \leq j'}^{a \equiv i} e_{[a,j')} e_{[i+i'-a;j)} - \sum_{i' \leq a < j'}^{a \equiv j} e_{[i;j+j'-a)} e_{[i',a)}\right)
\end{equation} 
$$
\frac {f_{[i;j)} f_{[i';j')}}{q^{\delta_{\barjj}^{\barii} - \delta_{\barii}^{\bari} - \delta_{\barjj}^{\bari}}} - \frac {f_{[i';j')} f_{[i;j)}}{q^{\delta_{\barjj}^{\barii} - \delta_{\barii}^{\barj} - \delta_{\barjj}^{\barj}}} = q \delta_{\bari}^{\barjj} f_{[i';j+j'-i)} - q \delta_{\barj}^{\barii} f_{[i+i'-j;j')} +
$$
\begin{equation}
\label{eqn:rtt 2}
+ (q-q^{-1}) \left(\sum_{i \leq a < j}^{a \equiv j'} f_{[i';j+j'-a)} f_{[i;a)} - \sum_{i<a\leq j}^{a \equiv i'} f_{[a;j)} f_{[i+i'-a;j')}\right) 
\end{equation}
as well as:
\begin{equation}
\label{eqn:rtt 3}
\Big[ e_{[i;j)}, f_{[i';j')} \Big] = \delta_{i-j}^{i'-j'} \delta_{\bari}^{\barii}  \cdot \frac {q^{h_i-h_j} - q^{h_j-h_i}}{q - q^{-1}} +
\end{equation}
$$
+ \begin{cases} \delta_{\barj}^{\barjj} \cdot e_{[i;j-j'+i')} \frac {q^{1 + h_{j'} - h_{i'}}}{q^{\delta_{\bari}^{\barii} - \delta_{\barii}^{\barj} + \delta_{\bari}^{\barj}}} - \delta_{\bari}^{\barii} \cdot e_{[i+j'-i';j)} q^{h_{i'} - h_{j'}} &\qquad \text{if } i-j<i'-j' \\ \delta_{\barj}^{\barjj} \cdot f_{[i';j'-j+i)} \frac {q^{h_{i} - h_{j}}}{q^{\delta_{\bari}^{\barii} - \delta_{\barii}^{\barj} + \delta_{\bari}^{\barj}}} - \delta_{\bari}^{\barii} \cdot f_{[i'+j-i;j')} q^{h_{j} - h_{i} - 1} & \qquad \text{if } i-j>i'-j' \end{cases} 
$$
$$
+ (q-q^{-1}) \sum_{k=1}^{\min(j-i,j'-i') - 1} \left( \delta_{\bark+\bari'}^{\barj} \frac {f_{[i'+k,j')} e_{[i;j-k)}}{q^{\delta_{\barii}^{\bari} - \delta_{\barii}^{\barj} + \delta_{\barj}^{\bari}}} \frac {\psi_{i'}}{\psi_{i'+k}} - \delta_{\bark+\bari}^{\barjj} \frac {e_{[i+k,j)} f_{[i';j'-k)}}{q^{\delta_{\bari}^{\barii} + \delta_{\bari}^{\barjj} - \delta_{\barii}^{\barjj}}} \frac {\psi_{i+k}}{\psi_{i}} \right) 
$$
We have rescaled the root generators \eqref{eqn:root} from the conventions in \cite{Tor} by a factor of $q^{-1} - q$ for the $e$'s and a factor of $1 - q^{-2}$ for the $f$'s. Therefore, we set:
\begin{equation}
\label{eqn:notation 0}
e_{[i;i)} = \frac 1{q^{-1} - q} \qquad \qquad f_{[i;i)} = \frac 1{1 - q^{-2}}
\end{equation}

\subsection{}\label{sub:verma} We have the triangular decomposition:
$$
\uu = \uug \otimes \uuo \otimes \uul
$$
where:
\begin{align*}
&\uug = \BC(q) \Big \langle e_{[i;j)} \Big \rangle_{i < j} \\ 
&\uuo = \BC(q) \Big \langle q^{h_k}, q^\gamma \Big \rangle_{1\leq k \leq n} \\ 
&\uul = \BC(q) \Big \langle f_{[i;j)} \Big \rangle_{i < j} 
\end{align*}
Moreover, $\uu$ is graded by $\nn$, where:
\begin{align*}
&\deg h_k = 0 \\
&\deg e_{[i;j)} = [i;j) \\
&\deg f_{[i;j)} = -[i;j)
\end{align*}
and we set $[i;j)  = \bs^i + ... + \bs^{j-1} \in \nn$. The algebra $\su$ has an analogous triangular decomposition, just by replacing the root generators $e_{[i;j)}$ (respectively the $f_{[i;j)}$) by the Drinfeld-Jimbo generators $e_i$ (respectively $f_i$) in all the formulas above. \\

\begin{definition}
\label{def:verma}

The {\bf universal Verma module}:
\begin{equation}
\label{eqn:verma}
M_{\hgl_n} = \uug \cdot v
\end{equation}
is the free $\uug$ module generated by a single vector $v$, which becomes a $\uu$ module if we impose additional relations:
\begin{align*}
f_{[i;j)} \cdot v &= 0 \\
q^{h_k} \cdot v &= u_k q^{-k} v \\
q^{\gamma} \cdot v &= q^n\oq v 
\end{align*}
$\forall i<j \in \BZ$ and $k \in \{1,...,n\}$, where $u_1,...,u_n$ and $\oq$ are formal parameters. \\

\end{definition}


\noindent Similarly, we may define universal Verma modules $M_{\hsl_n}$ and $M_{\hgl_1}$ by replacing: 
$$
\uu \text{ with } \su \text{ and } \uui, \text{ respectively}
$$
in Definition \ref{def:verma}. As a consequence of \eqref{eqn:quantum}, we have an isomorphism:
\begin{equation}
\label{eqn:verma decomp}
M_{\hgl_n} \cong M_{\hsl_n} \otimes M_{\hgl_1}
\end{equation}

\subsection{}
\label{sub:affdef}

We now recall the definition of affine Laumon spaces, for which we will use the notations in \cite{Aff}. Consider the surface $\BP^1 \times \BP^1$ and the divisors:
$$
D=\BP^1 \times \{0\}, \qquad \qquad \infty = \BP^1 \times \{\infty\} \cup \{\infty\} \times \BP^1.
$$
A \textbf{parabolic sheaf} $\CF_\bullet$ is a flag of rank $n$ torsion free sheaves: 
\begin{equation}
\label{eqn:flag}
\CF_n(-D) \subset \CF_1 ... \subset \CF_{n-1} \subset \CF_n
\end{equation}
on $\BP^1 \times \BP^1$, together with an isomorphism:
\begin{equation}
\label{eqn:framing}
\xymatrix{\CF_n(-D)|_\infty \ar[r] \ar[d]^\cong & \CF_{1}|_\infty \ar[r] \ar[d]^\cong & ... \ar[r] \ar[d]^\cong &  \CF_n|_\infty \ar[d]^\cong \\
\CO_\infty^{\oplus n}(-D) \ar[r] & \CO_\infty \oplus \CO_\infty^{\oplus n-1}(-D) \ar[r] & ... \ar[r]  & \CO_\infty^{\oplus n}}
\end{equation}
The data in \eqref{eqn:framing} is called \textbf{framing} at $\infty$, and it forces $c_1(\CF_i) = -(n-i)D$. On the other hand, $-c_2(\CF_i) =: d_i$ can vary over all non-negative integers, and we therefore call the vector $\bd = (d_1,...,d_{n}) \in \nn$ the \textbf{degree} of the parabolic sheaf $\CF_\bullet$. \\

\begin{definition}
\label{def:laumon}

The moduli space $\CM_{\bd}$ of rank $n$ degree $\bd$ parabolic sheaves is called an \textbf{affine Laumon space}. \\

\end{definition}

\noindent $\CM_\bd$ is a smooth quasiprojective variety of dimension $2|\bd| := 2(d_1+...+d_{n})$. The case of usual Laumon spaces is when $d_n = 0$, in which case parabolic sheaves \eqref{eqn:flag} reduce to a full flag of framed torsion-free sheaves on the divisor $D \cong \BP^1$. \\

\subsection{} 
\label{sub:equiv}

The maximal torus $T_n \subset GL_n$ and the rank 2 torus $\BC^* \times \BC^*$ act on $\CM_\bd$ by changing the trivialization at $\infty$, respectively by multiplying the base $\BP^1 \times \BP^1$. This allows us to define the $\tT-$equivariant $K-$theory group $K_{\tT}(\CM_{\bd})$, where $\tT$ is a $2^{n+2}-$fold cover of $T_n \times \BC^* \times \BC^*$. We will consider the representation ring of $\tT$:
\begin{equation}
\label{eqn:coordinates}
K_{\tT}(\pt) = \BC \left[u_1^{\pm 1},...,u_n^{\pm 1}, q^{\pm 1}, \oq^{\pm 1}\right] 
\end{equation}
The {\bf localized $K-$theory ring} is the $\text{Frac}(K_{\tT}(\pt))$--vector space:
$$
K_\bd := K_{\tT}(\CM_{\bd}) \bigotimes_{K_{\tT}(\pt)} \text{Frac}(K_{\tT}(\pt))
$$
and we will package these together, by setting:
\begin{equation}
\label{eqn:k}
K = \bigoplus_{\bd \in \nn} K_\bd
\end{equation}

\subsection{}\label{sub:corr} 

We will now introduce operators on $K$ which arise from geometric correspondences. For example, consider the {\bf fine correspondence} defined as:
\begin{multline}
\fZ_{[i;j)} := \Big\{ (\CF^j_\bullet \subset_x^{j-1} \CF^{j-1}_\bullet \subset_x^{j-2} ... \subset_x^{i+1} \CF^{i+1}_\bullet \subset^i_x \CF^i_\bullet), \text{ for some } x \in \BA^1 = D \backslash \infty \Big \} \\ \subset \CM_{\bd^j} \times \CM_{\bd^{j-1}} \times ... \times \CM_{\bd^{i+1}} \times \CM_{\bd^i}  \label{eqn:fine} 
\end{multline}
for arbitrary $\bd^+ = \bd^j, ..., \bd^i = \bd^- \in \nn$ such that $\bd^{k+1}  = \bd^k + \bs^k$. The notation:
$$
\CF'_\bullet \subset_x^k \CF_\bullet \quad \text{means that we have embeddings} \quad \{\CF'_i \hookrightarrow \CF_i\}_{1\leq i \leq n}
$$
which preserve the flag structure \eqref{eqn:flag} and the framing \eqref{eqn:framing}, such that all embeddings are isomorphisms, except for $\CF'_k \subset \CF_k$ which are such that $\CF_k/\CF'_k \cong \BC_x$. The fine correspondence is endowed with projection maps:
\begin{equation}
\label{eqn:fine maps}
\xymatrix{& \fZ_{[i;j)} \ar[ld]_{p^+} \ar[rd]^{p^-} & \\ \CM_{\bd^+} & & \CM_{\bd^-}}
\end{equation}
There exist line bundles $\CL_i,...,\CL_{j-1}$ on $\fZ_{[i;j)}$, whose fibers are given by:
\begin{equation}
\label{eqn:line bundles}
\CL_k|_{(\CF^j_\bullet \subset_x^{j-1} \CF^{j-1}_\bullet \subset_x^{j-2} ... \subset_x^{i+1} \CF^{i+1}_\bullet \subset^i_x \CF^i_\bullet)} = \Gamma(\BP^1 \times \BP^1, \CF^k_k/ \CF_k^{k+1})
\end{equation}
In \cite{Aff}, we also introduced the \textbf{eccentric correspondence} $\ofZ_{[i;j)}$, which does not admit a formula in terms of flags of sheaves as \eqref{eqn:fine}, but is endowed with maps:
\begin{equation}
\label{eqn:eccentric maps}
\xymatrix{& \ofZ_{[i;j)} \ar[ld]_{\op^+} \ar[rd]^{\op^-} & \\ \CM_{\bd^+} & & \CM_{\bd^-}}
\end{equation}
and with line bundles $\CL_1,...,\CL_{j-1}$ as in \eqref{eqn:line bundles}. Moreover, in \cite{Aff} we also introduced virtual structure sheaves (which should more appropriately be called derived scheme structures, but we will not go into these issues in the present paper):
\begin{align}
&[\fZ_{[i;j)}] \in K_T(\fZ_{[i;j)}) \label{eqn:virtual 1} \\
&[\ofZ_{[i;j)}] \in K_T(\ofZ_{[i;j)}) \label{eqn:virtual 2}
\end{align}
Let us consider the following modifications of these classes:
\begin{align}
&[\fZ_{[i;j)}^+] = [\fZ_{[i;j)}] \cdot u_{i+1}...u_j \cdot \frac {\CL_i}{\CL_{j-1}} \cdot (-1)^{j-i-1} \cdot \frac {q^{\left \lceil \frac {j-i}n \right \rceil-1}}{q^{d_i^+ - d_j^+}} \label{eqn:virtual mod 1} \\
&[\ofZ_{[i;j)}^-] = [\ofZ_{[i;j)}] \cdot \frac {u_{i}...u_{j-1}}{\CL_i ... \CL_{j-1}} \cdot (-1) \cdot \frac {q^{-j+i-\left \lfloor \frac {j-i}n \right \rfloor +2}}{q^{d_{i-1}^- - d_{j-1}^-}} \label{eqn:virtual mod 2}
\end{align}
Note that the right-hand sides of \eqref{eqn:virtual mod 1} and \eqref{eqn:virtual mod 2} differ by factors of $q$ and $-q^2$, respectively, from the same-named objects in \cite{Aff}. Moreover, the classes \eqref{eqn:virtual 1} and \eqref{eqn:virtual 2} differ by an additional factor of $1-q^2$ from the same-named objects in \cite{Aff}, because of the spaces $\fZ_{[i;j)}$ and $\ofZ_{[i;j)}$ studied in the present paper are equal to $\BA^1$ times the homonymous spaces in \cite{Aff} (while the point $x$ is allowed to range over $\BA^1$ in the current definition \eqref{eqn:fine}, we had fixed $x$ in \loccitt). \\

\subsection{}

Define the {\bf smooth correspondence}:
\begin{multline}
\fW_{\bd^+,\bd^-} := \Big\{ (\CF_\bullet^+ \subset \CF^-_\bullet) \text{ such that } \CF^-_i/\CF^+_i \text{ is scheme-theoretically} \\ \text{supported on } \{0\} \times \BP^1, \text{ for all } i \Big \}  \subset \CM_{\bd^+} \times \CM_{\bd^-}  \label{eqn:smooth}
\end{multline}
for all $\bd^\pm \in \nn$. We will write $\fW_k$ for the union on the correspondences \eqref{eqn:smooth} over all degree vectors such that $\bd^+ = \bd^- + k \bde$. There exists a rank $k$ vector bundle:
$$
\CL_i|_{(\CF^+_\bullet \subset \CF^-_\bullet)} = \Gamma(\BP^1 \times \BP^1, \CF_i^+/\CF_i^-)
$$ 
on $\fW_k$, for all $i \in \{1,...,n\}$. Consider the natural projection maps:
\begin{equation}
\label{eqn:smooth maps}
\xymatrix{& \fW_k \ar[ld]_{\pi^+} \ar[rd]^{\pi^-} & \\ \CM_{\bd^+} & & \CM_{\bd^-}}
\end{equation}
which remember $\CF_\bullet^+$ and $\CF_\bullet^-$, respectively. Schemes of the form \eqref{eqn:smooth} are well-known in geometric representation theory, as is the fact that they are smooth. However, in \cite{Aff}, we do not directly prove that \eqref{eqn:smooth} is actually smooth, but instead work with its virtual structure sheaf:
\begin{equation}
\label{eqn:virtual 3}
[\fW_k] \in K_T(\fW_k) 
\end{equation}
and the modifications:
\begin{align*}
&[\fW_k^+] = [\fW_k] \cdot \oq^{k^2-2k} \frac {(u_1u_2...u_{n})^k}{q^{nk}} \label{eqn:smooth mod 1} \\
&[\fW_k^-] = [\fW_k] \cdot \oq^{k^2} (-1)^{nk} \frac {(u_1u_2...u_{n})^k}{\prod_{i=1}^n \det \CL_i} \label{eqn:smooth mod 2}
\end{align*}
Consider the operators:
\begin{align}
&e_{[i;j)} : K \rightarrow K, \qquad \alpha \leadsto p^+_* \left( [\fZ_{[i;j)}^+] \cdot p^{-*} (\alpha) \right) \label{eqn:fine ops} \\
&f_{[i;j)} : K \rightarrow K, \qquad \alpha \leadsto \op^-_* \left( [\ofZ_{[i;j)}^-] \cdot \op^{+*} (\alpha) \right) \label{eqn:eccentric ops} \\
&g_{\pm k} \ : K \rightarrow K, \qquad \alpha \leadsto \pi^\pm_* \left( [\fW_k^\pm] \cdot \pi^{\mp*} (\alpha) \right) \label{eqn:smooth ops}
\end{align}
The main result of \cite{Aff} is the following: \\

\begin{theorem}
\label{thm:geom} 

There is an action $\uu \curvearrowright K$, where: \\

\begin{itemize}[leftmargin=*]
	
\item the elements $e_{[i;j)}$, $f_{[i;j)}$ of \eqref{eqn:root} act by the operators \eqref{eqn:fine ops}, \eqref{eqn:eccentric ops}, respectively \\

\item the central element acts by $q^\gamma = q^n\oq$, and:
\begin{equation}
\label{eqn:psi}
q^{h_i} = \text{multiplication by } u_i q^{d _i - d_{i-1} - i} \text{ on } K_\bd
\end{equation}

\item the elements $g_{\pm k} \in \uui \subset \uu$ of \eqref{eqn:group-like} act on $K$ by the operators \eqref{eqn:smooth ops} \\

\end{itemize}

\noindent Moreover, we have an isomorphism of $\uu$--modules $K \cong M_{\hgl_n}$. \\	

\end{theorem} 

\begin{remark}
\label{rem:fix constant}
	
For the Theorem to hold as stated, the generating series \eqref{eqn:group-like}	of the elements $g_{\pm k}$ must be defined as the exponential of the generating series of elements $p_{\pm k}$ which satisfy the following relation:
\begin{equation}
\label{eqn:fix constant}
[p_k, p_l] = \delta_{k+l}^0 k \cdot \frac {q^{nk} - q^{-nk}}{\oq^k - \oq^{-k}}
\end{equation}
In other words, formula \eqref{eqn:fix constant} and $q^\gamma = q^n \oq$ fix the constant in \eqref{eqn:heis rel}. \\
	
\end{remark}

\subsection{}\label{sub:e}

Let us consider two different copies of affine Laumon space, henceforth denoted by $\CM_{\bd}$ and $\CM_{\bd'}$, which can carry different sets of equivariant parameters $u_1,...,u_n$ and $u_1',...,u_n'$, respectively. More rigorously, we consider the action of the bigger torus $(\BC^*)^n \times (\BC^*)^n \times \BC^* \times \BC^*$ on the product $\CM_{\bd} \times \CM_{\bd'}$, and write:
\begin{align*}
&K = \bigoplus_{\bd \in \nn} K_{\bd} \quad \text{as a module over} \ \ \quad \BC[u_1^{\pm 1},...,u_n^{\pm 1}, q, \oq] \\
&K' = \bigoplus_{\bd' \in \nn} K_{\bd'} \quad \text{as a module over} \quad \BC[{u_1'}^{\pm 1},...,{u_n'}^{\pm 1}, q, \oq]
\end{align*}
Consider the following vector bundle on $\CM_{\bd} \times \CM_{\bd'}$:
\begin{equation}
\label{eqn:e}
\CE|_{(\CF_\bullet,\CF_\bullet')} = \textrm{Ext}^1(\CF'_\bullet,\CF_\bullet(-\infty)).
\end{equation}
The actual definition of $\textrm{Ext}$ of two flags of sheaves is given in \cite{FFNR}, where it is shown that \eqref{eqn:e} is a vector bundle. We may consider the full exterior power of \eqref{eqn:e}:
$$
[\wedge^\bullet (\CE^\vee,m)] = \sum_{k = 0}^{\text{rank }\CE} \left( - q^{2m} \right)^k [\wedge^k \CE^{\vee}] \in K_T(\CM_{\bd} \times \CM_{\bd'})
$$ 
Consider also the following modification of the class above:
\begin{multline}
[\widetilde{\wedge}^\bullet (\CE^\vee,m)] = [\wedge^\bullet (\CE^\vee,m)] \cdot \\ \cdot (-1)^{|\bd|} \frac {u_{\bd'+1}}{u_{\bd}} q^{-\langle \bd',\bd \rangle - |\bd'|-m (|\bd|+|\bd'|)} \prod_{i=1}^n \det \CV_i \cdot \lambda_{\bd'} \label{eqn:modification} 
\end{multline}
where the tautological vector bundle $\CV_i$ has rank $d_i$ on $\CM_{\bd} \subset \CM$, and:
\begin{equation}
\label{eqn:lambda}
\lambda_{\bd'} = (-1)^{|\bd'|} \frac {u_{\bd'}}{u_{\bd'+1}} \frac {q^{|\bd'| + \langle \bd',\bd' \rangle}}{\det \CV_i'}
\end{equation} 
Then let us define the operator:
\begin{equation}
\label{eqn:op}
A_m : K_{\bd'} \longrightarrow K_{\bd}
\end{equation}
$$
\alpha \leadsto p_{1*} \Big( [\widetilde{\wedge}^\bullet(\CE^\vee,m)] \cdot p^{*}_2(\alpha)\Big)
$$
where $p_1,p_2$ are the projections from $\CM_{\bd} \times \CM_{\bd'}$ onto the first and second factors, respectively. If we consider all $\bd,\bd' \in \nn$ together, we obtain an operator: 
\begin{equation}
\label{eqn:A}
A_m : K' \longrightarrow K
\end{equation}
Strictly speaking, the operator \eqref{eqn:A} maps into the completion of $K = \oplus_{\bd \in \nn} K_{\bd}$. \\

\subsection{}\label{sub:commutation} 

The main application of representation theory to the study of the Ext operator $A_m$ is to connect the latter to intertwiners for quantum affine algebras. Our main geometric computations, to be proved in Section \ref{sec:inter}, are given below: \\

\begin{proposition}
\label{prop:comm root e}
The operator $A_m$ satisfies the following commutation relations:
\begin{multline}
e_{[i;j)} A_m + e_{[i+1;j)} q^{-h_{i+1}} A_m q^{h_i+m} = A_m e_{[i;j)} + q^{h_j} A_m q^{-h_{j-1}-m} e_{[i;j-1)} \label{eqn:comm1}
\end{multline}
for all $i<j$. \\ 
\end{proposition}

\noindent In particular, since $e_{[i;i)} = \frac 1{q^{-1}-q}$, we have the following relation in the case $j=i+1$:
\begin{equation}
\label{eqn:special 1}
e_i A_m - A_m e_i = \frac {q^{-h_{i+1}} A_m q^{h_i + m} - q^{h_{i+1}} A_m q^{-h_i - m}}{q-q^{-1}}
\end{equation}
$$$$

\begin{proposition}
\label{prop:comm root f}
The operator $A_m$ satisfies the following commutation relations:
\begin{multline}
f_{[i+1,j+1)} A_m + \sum_{a=i+1}^j (-q)^{i-a} r_{\left \lfloor \frac {a-i}n \right \rfloor} f_{[a+1;j+1)} q^{-h_{i+1}-...-h_a} A_m q^{h_{i+1}+...+h_a+m(a-i)} = \\ = A_m f_{[i;j)} + \sum_{a=i}^{j-1} (-q)^{a-j} r_{\left \lfloor \frac {j-a}n \right \rfloor} q^{h_{a+1}+...+h_j} A_m q^{-h_{a+1}-...-h_j-m(j-a)} f_{[i;a)} \label{eqn:comm2}
\end{multline}
for all $j>i$, where we set $r_k = \prod_{l=1}^k \frac {q \oq^{2l} - q^{-1}}{\oq^{2l}-1}$. \\

\end{proposition}

\noindent In particular, since $f_{[i;i)} = \frac 1{1 - q^{-2}}$, we have the following relation in the case $j=i+1$:
\begin{equation}
\label{eqn:special 2}
f_i A_m - A_m f_{i-1} = \frac {q^{-h_i} A_m q^{h_i + m} - q^{h_i} A_m q^{-h_i-m}}{q-q^{-1}}
\end{equation}
$$$$

\begin{proposition}
\label{prop:comm group}
	
The operator $A_m$ satisfies the following commutation relations:
\begin{align}
&[p_k, A_m] \ = (-1)^{k-1} \cdot \frac {q^{(m+1)nk} \oq^k - q^{-(m+1)nk} \oq^{-k}}{\oq^{k}-\oq^{-k}} \cdot A_m \label{eqn:comm3} \\
&[p_{-k}, A_m] = (-1)^{k-1} \cdot \frac {q^{mnk} - q^{-mnk}}{\oq^{k}-\oq^{-k}} \cdot A_m \label{eqn:comm4}
\end{align}
for all $k > 0$. \\

\end{proposition}

\subsection{}

The decomposition \eqref{eqn:verma decomp} and the last sentence of Theorem \ref{thm:geom} imply that:
$$
K \cong M_{\hgl_n} = M_{\hsl_n} \otimes M_{\hgl_1}
$$
and the action of $\uu$ on the LHS is matched with the action of $\su \otimes \uui$ on the RHS. Because relations \eqref{eqn:special 1}, \eqref{eqn:special 2} do not involve the $\uui$ generators $p_{\pm k}$, and relations \eqref{eqn:comm3}, \eqref{eqn:comm4} do not involve the $\su$ generators $e_i, f_i$, then:
\begin{equation}
\label{eqn:tensor}
A_m = A_m^{\hsl_n} \otimes A_m^{\hgl_1} :  M' \rightarrow M
\end{equation}
where:
\begin{equation}
\label{eqn:prime}
M_{\hsl_n}' \stackrel{A_m^{\hsl_n}}\longrightarrow M_{\hsl_n} \quad \text{satisfies relations \eqref{eqn:special 1}, \eqref{eqn:special 2}}
\end{equation}
\begin{equation}
\label{eqn:double prime}
M_{\hgl_1}' \stackrel{A_m^{\hgl_1}}\longrightarrow M_{\hgl_1} \quad \text{satisfies relations \eqref{eqn:comm3}, \eqref{eqn:comm4}}
\end{equation}
(for any symbol $* \in \{\hgl_n, \hsl_n, \hgl_1\}$, we write $M_*$ and $M_*'$ for two copies of the Verma module, with potentially different generic highest weights $u_1,...,u_r$ and $u_1',...,u_r'$). The following is a well-known property of Verma modules. \\

\begin{proposition}
\label{prop:unique}

There is a unique (up to scalar) operator $M'_* \stackrel{A_m}\longrightarrow M_*$ satisfying:
$$
\begin{cases} 
\text{\eqref{eqn:comm1} and \eqref{eqn:comm2}} &\text{if } * = \hgl_n \\
\text{\eqref{eqn:special 1} and \eqref{eqn:special 2}} &\text{if } * = \hsl_n \\
\text{\eqref{eqn:comm3} and \eqref{eqn:comm4}} &\text{if } * = \hgl_1
\end{cases}
$$

\end{proposition}

\begin{proof} We will only prove the $* = \hgl_n$ case, and so we abbreviate $M = M_{\hgl_n}$ from now on (we leave the remaining two cases as exercises for the interested reader). Proposition 2.17 of \cite{Tor} establishes the fact that a linear basis of the universal Verma module $M$ (over the ground field $\BQ(u_1,...,u_n,q,\oq)$) consists of expressions:
\begin{equation}
\label{eqn:basis}
e_{[i_1;j_1)}... e_{[i_d;j_d)} \cdot v
\end{equation}
over all $\{i_s < j_s\}_{1\leq s \leq d} \subset \{1,...,n\} \times \BZ$, ordered such that $j_1-i_1 \leq ... \leq j_d-i_d$ (and if $j_s-i_s = j_{s+1}-i_{s+1}$, then we require $i_s \leq i_{s+1}$). Iterating \eqref{eqn:comm1} yields:
$$
A_m e_{[i;j)} = \sum_{\e \in \{0,1\}} \sum_{a=i}^{j-\e} ... e_{[a;j-\e)} A_m ...
$$
where $...$ stand for various constants and products of Cartan elements $q^{\pm h_k}$. With \eqref{eqn:basis} in mind, we conclude that knowing the vector $v_m = A_m \cdot v$ completely determines the operator $A_m$. However, applying \eqref{eqn:comm2} to the vector $v$ implies:
$$
f_{[i+1;j+1)} \cdot v_m = \sum_{a=i+1}^j ... f_{[a+1;j+1)} \cdot v_m
$$
where $...$ stand for various constants and products of Cartan elements $q^{\pm h_k}$. Iterating the relation above implies that: 
\begin{equation}
\label{eqn:iteration}
f_{[i;j)} \cdot v_m = P_{ij}(q^{h_1},...,q^{h_n})\cdot v_m
\end{equation}
for some Laurent polynomials $P_{ij}$ that one can deduce algorithmically. Hence, if:
$$
v_m = \alpha v + \sum_{d=1}^\infty \sum^{\{i_s < j_s\}}_{\text{as in \eqref{eqn:basis}}} \alpha_{i_1,j_1,...,i_d,j_d} e_{[i_1;j_1)}... e_{[i_d;j_d)} \cdot v
$$
then the constants $\alpha_{i_1,j_1,...,i_d,j_d}$ can be computed from $\alpha$ by induction on $j_i-i_1+...+j_d-i_d$ by using \eqref{eqn:rtt 3}, \eqref{eqn:iteration}, and the fact that there exists a non-degenerate bialgebra pairing between the span of ordered products of $e_{[i;j)}$'s and the span of ordered products of $f_{[i;j)}$'s (see \cite{Tor}). This implies that $v_m$ is unique up to scalar. 
	
\end{proof}

\subsection{} 

We expect that the operator $M' \stackrel{A_m}\longrightarrow M$ defined by \eqref{eqn:comm1}, \eqref{eqn:comm2} can be related to intertwiners between Verma modules of $\uu$, and in the next Section we will prove this fact in the limit when the quantum group $\uu$ is replaced by the affine Lie algebra $\hgl_n$. In the meantime, we observe that \eqref{eqn:tensor} reduces the problem from $\uu$ to $\su$, because it is easy to show that:
\begin{multline}
A_m^{\hgl_1} = \exp \left( \sum_{k=1}^\infty \frac {p_k (-1)^k}k \cdot \frac {q^{mnk} - q^{-mnk}}{q^{nk} - q^{-nk}} \right) \cdot \\ \cdot \exp \left(- \sum_{k=1}^\infty \frac {(-1)^{k} p_{-k}}k \cdot \frac {q^{(m+1)nk} \oq^k - q^{-(m+1)nk} \oq^{-k}}{q^{nk} - q^{-nk}} \right) 
\label{eqn:vertex}
\end{multline}
Indeed, by the uniqueness statement of Proposition \ref{prop:unique}, it is enough to show that the two sides of \eqref{eqn:vertex} satisfy the same commutation relations with the generators $p_{\pm k}$. In the case of the LHS, these are given by \eqref{eqn:comm3} and \eqref{eqn:comm4}. Meanwhile:
\begin{align*}
&[p_k, \text{RHS of \eqref{eqn:vertex}}] = (-1)^{k-1} \cdot \frac {q^{(m+1)nk} \oq^k - q^{-(m+1)nk} \oq^{-k}}{\oq^{k}-\oq^{-k}} \cdot \text{RHS of \eqref{eqn:vertex}} \\
&[p_{-k}, \text{RHS of \eqref{eqn:vertex}} ] = (-1)^{k-1} \cdot \frac {q^{mnk} - q^{-mnk}}{\oq^{k}-\oq^{-k}} \cdot \text{RHS of \eqref{eqn:vertex}}
\end{align*}
are straightforward consequences of \eqref{eqn:fix constant}. \\

\begin{remark}

Note that, a priori, \eqref{eqn:vertex} only holds up to scalar multiple. However, in \eqref{eqn:tensor} we can freely multiply $A_m^{\hgl_1}$ by any constant, at the cost of multiplying $A_m^{\hsl_n}$ by the inverse constant. Therefore, we may take \eqref{eqn:vertex} as a definition of $A_m^{\hgl_1}$. \\

\end{remark}

\subsection{}

When $u_i = u_i'$, $\forall i \in \{1,...,n\}$, the groups $K \cong M$ and $K' \cong M'$ are naturally identified, and we can talk about the graded character of the operator $A_m$:
$$
\chi_{A_m} = \sum_{\bd \in \nn} \bx^{\bd} \cdot \text{Tr} \Big( \textrm{proj}_{K_\bd}(A_m|_{K_\bd}) \Big)
$$
where $\bx^{\bd} = x_1^{d_1}... x_n^{d_n}$ are formal variables. \\

\begin{proposition} \label{prop:char}
	
We have $\chi_{A_m} = Z_m$. \\
	
\end{proposition}

\begin{proof} The proof closely follows that of Proposition 6.6 of \cite{Laumon}, so we leave the details to the interested reader. The main thing it uses is that the restriction of the Ext bundle to the diagonal $\Delta \cong \CM_{\bd} \hookrightarrow \CM_{\bd} \times \CM_{\bd}$ is given by:
$$
\widetilde{\wedge}^\bullet(\CE,m) \Big|_\Delta = \wedge^\bullet(\CE,m) \Big|_\Delta q^{-2m|\bd|} = \wedge^{\bullet}(\text{Tan}^\vee \CM_{\bd},m) q^{-2m|\bd|}
$$
The second equality can be seen, for example, by setting $\CV_i  = \CV_i'$ in \eqref{eqn:ke}.

	
\end{proof}

\noindent Therefore, computing the partition function of 5d ${\mathcal{N}}=1$ supersymmetric gauge theory with adjoint matter (in the presence of a surface operator insertion) amounts to computing the graded character of the operator $A_m : M \rightarrow M$, which is completely determined by properties \eqref{eqn:comm1} and \eqref{eqn:comm2}. By \eqref{eqn:tensor}, we have :
$$
\chi_{A_m} = \chi_{A^{\hsl_n}_m} \chi_{A^{\hgl_1}_m}
$$
To compute $\chi_{A^{\hsl_n}_m}$ from \eqref{eqn:special 1}, \eqref{eqn:special 2} would require one to relate $A^{\hsl_n}_m$ to $\su$ intertwiners, and then to apply \cite{ESV}, \cite{S} in order to compute characters of intertwiners. \\

\begin{proposition}
\label{prop:char vertex}

The character of the operator \eqref{eqn:vertex} is given by:
\begin{equation}
\label{eqn:char vertex}
\chi_{A^{\hgl_1}_m} = \exp \left[\sum_{k=1}^\infty \frac{x^k}{(1-x^k)k} \left(1 + \frac {(q^{(m+1)nk}\oq^k - q^{-(m+1)nk}\oq^{-k})(q^{mnk} - q^{-mnk})}{(\oq^k - \oq^{-k})(q^{nk} - q^{-nk})} \right) \right] \qquad
\end{equation}

\end{proposition}

\begin{proof} We will denote partitions by $\lambda = (1^{n_1}2^{n_2}...)$ and $\mu = (1^{m_1}2^{m_2}...)$. Define:
$$
|\lambda| = n_1+2n_2+..., \qquad z_\lambda = \prod_{k=1}^\infty k^{n_k} n_k!, \qquad {\lambda \choose \mu} = \prod_{k=1}^\infty {n_k \choose m_k}
$$
and write $p_{\pm \lambda} = p_{\pm 1}^{n_1} p_{\pm 2}^{n_2} ...$. The Verma module of the Heisenberg algebra is called the Fock space, and we will denote it by $\Lambda$. It has a basis consisting of monomials:
\begin{equation}
\label{eqn:basis heis}
p_\lambda \cdot v = p_1^{n_1}p_2^{n_2}... \cdot v
\end{equation}
as $\lambda$ goes over all partitions. We define $\Lambda_d \subset \Lambda$ to be the linear span of the vectors $p_\lambda \cdot v$ with $|\lambda|=d$. To keep the notation simple, let us rescale the generators of the Heisenberg algebra so that they satisfy the commutation relations $[p_k,p_l] = -\delta_{k+l}k$. The problem then becomes to show that:
$$
B = \exp \left(\sum_{k=1}^\infty \frac {\alpha_k p_k}k \right) \exp \left(\sum_{k=1}^\infty \frac {\beta_k p_{-k}}k \right) = \sum^{\mu,\nu}_{\text{partitions}} \frac {\alpha_\mu \beta_\nu}{z_\mu z_\nu} p_\mu p_{-\nu}
$$
(where $\alpha_\lambda = \prod_{k=1}^\infty \alpha_k^{n_k}$ for $\lambda = (1^{n_1}2^{n_2}...)$) has graded character:
\begin{equation}
\label{eqn:char heis}
\chi_B := \sum_{d=0}^\infty x^d \Tr \left( B|_{\Lambda^d} \right) = \exp \left[\sum_{k=1}^\infty \frac {x^k  (1+\alpha_k\beta_k)}{k(1-x^k)} \right]
\end{equation}
One can prove the equality above by showing that the character of $B$ satisfies the differential equation:
$$
\left[ x \frac {\partial}{\partial x} -  \sum_{k=1}^\infty \frac {x^k(1+\alpha_k\beta_k)}{(1-x^k)^2} \right] \chi_B = 0
$$
which is done along the lines of Proposition \ref{prop:char inter}. However, we can also prove formula \eqref{eqn:char heis} by explicitly computing the trace in the basis \eqref{eqn:basis heis}. It is a straightforward exercise, which we leave to the interested reader, to show that:
$$
p_{-\nu} p_\lambda \cdot v = \begin{cases} {\lambda \choose \nu} z_\nu p_{\lambda \backslash \nu} \cdot v &\text{if } \nu \subset \lambda \\ 0 &\text{otherwise } \end{cases}
$$
where we write $\nu \subset \lambda$ if all the parts of the partition $\nu$ (with multiplicities) are also in the partition $\lambda$. Clearly, $p_\mu p_{\lambda \backslash \nu}$ is a multiple of $p_\lambda$ only if $\mu = \nu$, hence:
$$
\chi_B = \sum^{\lambda}_{\text{partition}} x^{|\lambda|} \sum_{\mu \subset \lambda} \frac {\alpha_\mu \beta_\mu}{z_\mu} {\lambda \choose \mu} = \sum_{\mu = (1^{m_1}2^{m_2}...)} x^{|\mu|} \frac {\alpha_\mu \beta_\mu}{z_\mu} \prod_{k=1}^\infty \sum_{d_k = 0}^\infty x^{kd_k} {m_k + d_k \choose m_k} 
$$
$$
= \sum_{\mu = (1^{m_1}2^{m_2}...)} x^{|\mu|} \frac {\alpha_\mu \beta_\mu}{z_\mu} \prod_{k=1}^\infty (1-x^k)^{-m_k-1} = \prod_{k=1}^\infty \frac 1{1-x^k} \left(\sum_{m=0}^\infty \frac {x^{mk}}{(1-x^k)^m} \cdot \frac {(\alpha_k\beta_k)^m}{k^m m!} \right)
$$
The right-hand side above is easily seen to match the right-hand side of \eqref{eqn:char heis}. 

\end{proof}

\section{Cohomology and the Main Theorem}
\label{sec:proof}


\subsection{}\label{sub:lie1}

In the limit $q\rightarrow 1$, the quantum group $\su$ degenerates to the universal enveloping of the {\bf Kac-Moody Lie algebra}:
$$
U(\hsl_n) := \BC \Big \langle \oe_1, ..., \oe_n, \of_1, ..., \of_n, h_1...,h_{n} \Big \rangle = \lim_{q\rightarrow 1} \su
$$ 
where $\oe_i = \lim_{q \rightarrow 1} e_i$, $\of_i = \lim_{q \rightarrow 1} f_i$. The commutation relations between the generators above are obtained by taking the $q \rightarrow 1$ limit of relations \eqref{eqn:sug1}--\eqref{eqn:sug4}:
\begin{align}
	&[h_j, \oe_i] = + \langle \bs^i, \bs^j \rangle \oe_i \label{eqn:sug1 lie} \\
	&[h_j, \of_i] = - \langle \bs^i, \bs^j \rangle \of_i \label{eqn:sug1 prim lie}
\end{align}
\begin{equation}
\label{eqn:sug2 lie}
\left [\oe_i, \of_j \right ] = \delta_j^i (h_i - h_{i+1})
\end{equation}
as well as:
\begin{equation}
\label{eqn:sug3 lie}
[\oe_i, \oe_j] = 0, \ \text{if } i - j \not \equiv \{-1,1\} \text{ mod } n
\end{equation}
\begin{equation}
\label{eqn:sug4 lie}
\left[ \oe_i, [\oe_i, \oe_{i + 1}] \right] = \left[ \oe_i, [\oe_i, \oe_{i-1}] \right] = 0
\end{equation}
and the analogous relations for $\of$'s instead of $\oe$'s. Note that in the limit $q\rightarrow 1$, the coproduct \eqref{eqn:cop1}--\eqref{eqn:cop3} tends to the standard Lie algebra coproduct:
$$
\Delta(x) = x\otimes 1 + 1 \otimes x, \qquad \forall x\in \hsl_n
$$

\subsection{}
\label{sub:above}

Similarly, we may consider the generators \eqref{eqn:root} of $\uu$, and define:
$$
\oe_{[i;j)} = \lim_{q\rightarrow 1} e_{[i;j)} \qquad \qquad \of_{[i;j)} = \lim_{q\rightarrow 1} f_{[i;j)} 
$$
which satisfy the following $q \rightarrow 1$ limits of relations \eqref{eqn:rtt 1}--\eqref{eqn:rtt 3}:
\begin{equation}
\label{eqn:cohom comm 1}
\Big[ \oe_{[i;j)}, \oe_{[i';j')} \Big] =\delta_{\barii}^{\barj} \oe_{[i;j+j'-i')} - \delta_{\bari}^{\barjj} \oe_{[i+i'-j';j)}
\end{equation}
\begin{equation}
\label{eqn:cohom comm 2}
\Big[ \of_{[i;j)}, \of_{[i';j')} \Big] = \delta_{\bari}^{\barjj} \of_{[i';j+j'-i)} - \delta_{\barj}^{\barii} \of_{[i+i'-j;j')}
\end{equation}
and:
\begin{equation}
\label{eqn:cohom comm 3}
\Big[ \oe_{[i;j)},  \of_{[i';j')}\Big] = \delta_{i-j}^{i'-j'} \delta_{\bari}^{\barii} \cdot (h_i-h_j) +
\end{equation}
$$
+ \begin{cases} \delta_{\barj}^{\barjj} \cdot \oe_{[i;j-j'+i')} - \delta_{\bari}^{\barii} \cdot \oe_{[i+j'-i';j)} & \qquad \text{if } i - j < i' - j' \\ \delta_{\barj}^{\barjj} \cdot \of_{[i';j'-j+i)} - \delta_{\bari}^{\barii} \cdot \of_{[i'+j-i;j')} & \qquad \text{if } i-j>i'-j' \end{cases}
$$
We will henceforth write $M$ for the universal Verma module, which is freely generated by the subalgebra $\langle \oe_{[i;j)} \rangle_{i<j}$ acting on a vector $v$, subject to the relations:
\begin{align*}
\of_{[i;j)} \cdot v &= 0 \\ h_k \cdot v &= (a_k - k)v \\ \gamma \cdot v &= (n+p)v 
\end{align*}
where $a_1,...,a_n, p$ are formal parameters (specifically, $a_i = \log_q u_i$ and $p = \log_q \oq$). \\


\subsection{}

Consider the affine Lie algebra:
$$
\hgl_n = \Big(n \times n \text{ matrices valued in }\BC[z,z^{-1}] \Big) \oplus \BC \cdot \gamma
$$
endowed with the usual Lie algebra bracket that has $\gamma$ central, and:
\begin{equation}
\label{eqn:lie affine}
[X \cdot z^k, Y \cdot z^l] = [X,Y] \cdot z^{k+l} + \delta_{k+l}^0 k \cdot \Tr(XY)\gamma
\end{equation}
Let $E_{ij}$ denote the elementary matrix with a single 1 at the intersection of row $i$ and column $j$, and 0 everywhere else. We introduce the notation:
\begin{equation}
\label{eqn:notation}
E_{ij} = E_{\bari \barj} \cdot z^{\left \lfloor \frac {j-1}n \right \rfloor - \left \lfloor \frac {i-1}n \right \rfloor}
\end{equation} 
for all integers $i,j$. With this in mind, it is easy to see that the assignment:
\begin{align}
\oe_{[i;j)} &\leadsto E_{ij} \label{eqn:ass 1} \\
\of_{[i;j)} &\leadsto E_{ji} \label{eqn:ass 2} \\
h_i &\leadsto E_{ii} - \gamma \left \lfloor \frac {i-1}n \right \rfloor \label{eqn:ass 3} 
\end{align}
$\forall i<j \in \BZ$, matches relations \eqref{eqn:cohom comm 1}--\eqref{eqn:cohom comm 3} with the Lie bracket \eqref{eqn:lie affine}. Moreover:
\begin{equation}
\label{eqn:power}
P_k = \sum_{i=1}^n E_{i,i+nk} = \text{identity matrix} \cdot z^k  
\end{equation}
$\forall k \in \BZ \backslash 0$, satisfy the relations:
\begin{equation}
\label{eqn:heis lie}
[P_k, P_l] = \delta_{k+l}^0 kn\gamma 
\end{equation}
and are thus limits $q\rightarrow 1$ of the elemnts $p_{\pm k}$ of Subsection \ref{sub:quantum2}, properly rescaled. \\

\subsection{} Let us now consider the degeneration of equivariant $K$--theory to equivariant cohomology. Specifically, rename the equivariant parameters of Subsection \ref{sub:equiv} to:
$$
q = e^{\hbar}, \quad \oq = e^{\hbar p}, \quad u_i = e^{\hbar a_i} \ \forall i \in \{1,...,n\}
$$
and let us consider only the leading order terms in $\hbar$ in all our formulas. Expressions such as $\prod (1-q^{-i}\oq^{-j}u_1^{-k_1}... u_n^{-k_n})$ have leading order term $\prod (i + jp + a_1k_1 + ... + a_n k_n)$ as $\hbar \rightarrow 0$. Moreover, if we define the equivariant cohomology groups:
$$
H = \bigoplus_{\bd \in \nn} H_{\bd} \qquad \text{where} \qquad H_{\bd} = H_T(\CM_{\bd}) \bigotimes_{H_T(\pt)} \text{Frac}(H_T(\pt))
$$
then the operators \eqref{eqn:fine ops} and \eqref{eqn:eccentric ops} degenerate to operators:
\begin{align*}
&\oe_{[i;j)} : H \rightarrow H, \qquad \alpha \leadsto p^+_* \left( [\fZ_{[i;j)}] \cdot p^{-*} (\alpha) \right)  \cdot (-1)^{j-i-1}  \\
&\of_{[i;j)} : H \rightarrow H, \qquad \alpha \leadsto \op^-_* \left( [\ofZ_{[i;j)}] \cdot \op^{+*} (\alpha) \right) \cdot (-1)
\end{align*}
which satisfy relations \eqref{eqn:cohom comm 1}--\eqref{eqn:cohom comm 3}, and thus induce an isomorphism:
\begin{equation}
\label{eqn:iso cohom}
H \cong M
\end{equation}
with the universal Verma module defined at the end of Subsection \ref{sub:above}. Similarly, the operator \eqref{eqn:A} degenerates to an operator:
\begin{equation}
\label{eqn:oper cohom}
\tA_m : H' \rightarrow H
\end{equation}
where $H$ and $H'$ denote two copies of the cohomology groups of affine Laumon spaces, with equivariant parameters $a_1,...,a_n$ and $a_1',...,a_n'$, respectively. \\ 

\subsection{}

By taking the limit $q \rightarrow 1$ in Propositions \ref{prop:comm root e} and \ref{prop:comm root f}, we note that the operator \eqref{eqn:oper cohom} satisfies the following commutation relations with the root generators:
\begin{equation}
\label{eqn:cohom comm1} 
\left( \oe_{[i;j)} + \oe_{[i+1;j)} \right) \tA_m = \tA_m \left( \oe_{[i;j)} + \oe_{[i;j-1)}  + m \delta_j^{i+1} \right)
\end{equation}
(we will use the notation $\oe_{[i;i)} = \of_{[i;i)} = h_i$ in the present Subsection, in order to keep all our formulas concise) for all $i<j$, and:
\begin{multline}
\left[ \sum_{a=i}^{j-1} (-1)^{i-a} \orr_{\left \lfloor \frac {a-i}n \right \rfloor} \of_{[a+1;j+1)} - (-1)^{i-j} \orr_{\left \lfloor \frac {j-i}n \right \rfloor} (h_{i+1}+...+h_j) \right] \tA_m = \\ = \tA_m \left[ \sum_{a=i+1}^{j} (-1)^{a-j} \orr_{\left \lfloor \frac {j-a}n \right \rfloor} \of_{[i;a)} - (-1)^{i-j} \orr_{\left \lfloor \frac {j-i}n \right \rfloor}(h_{i+1}+...+h_j + m (j-i)) \right] \qquad \label{eqn:cohom comm2}
\end{multline}
for all $i<j$, where $\orr_k = \prod_{l=1}^k \frac {1+pl}{pl} = {k + \frac 1p \choose k }$. \\

\begin{proposition}
\label{prop:matrix notation}

In the matrix notation \eqref{eqn:ass 1}--\eqref{eqn:ass 3}, the relations above read:
\begin{equation}
\label{eqn:cohom}
\left( E_{ij} + E_{i+1,j} \right) \tA_m = \tA_m \left( E_{ij} + E_{i,j-1} + \alpha_{i-j}m \right)
\end{equation}
for all $i,j \in \BZ$, where:
\begin{equation}
\label{eqn:alpha}
\alpha_k = \begin{cases} 1 & \text{if } k \in \{0,-1\}  \\  (-1)^{ns-1} \frac np   & \text{if } k = ns \text{ or } ns-1 \text{ with } s>0 \\ 0 & \text{for all other integers} \end{cases}
\end{equation}

\end{proposition}


\begin{proof} It is clear that \eqref{eqn:cohom comm1} is precisely \eqref{eqn:cohom} when $i<j$. In order to prove the case $i \geq j$, let us add together formula \eqref{eqn:cohom comm2} for the pairs $(i,j)$ and $(i,j-1)$:
\begin{equation}
\label{eqn:two line}
\left[ \sum_{a = i}^{j-1} (-1)^{i-a} \orr_{\left \lfloor \frac {a-i}n \right \rfloor} \left(\of_{[a+1;j+1)} + \of_{[a+1;j)}\right) \right] \tA_m = 
\end{equation}
$$
= \tA_m \left[ \of_{[i;j)} + \sum_{i < a < j}^{n|j-a}  \frac {(-1)^{a-j}  \orr_{\left \lfloor \frac {j-1-a}n \right \rfloor} \of_{[i;a)} }{p\left \lfloor \frac {j-a}n \right \rfloor} \right. + 
$$
$$
\left. + (-1)^{i-j+1} \left( \orr_{\left \lfloor \frac {j-i-1}n \right \rfloor} h_j + \orr_{\left \lfloor \frac {j-i}n \right \rfloor} m(j-i) - \orr_{\left \lfloor \frac {j-i-1}n \right \rfloor} m(j-i-1) \right) \right]
$$
(above, we used the simple identity $\orr_{k} - \orr_{k-1} = \frac {\orr_{k-1}}{pk}$, as well as the formula:
$$
\left( \orr_{\left \lfloor \frac {j-i}n \right \rfloor} - \orr_{\left \lfloor \frac {j-1-i}n \right \rfloor} \right)[h_{i+1}+...+h_j, \tA_m] = 0
$$
which happens because $h_{i+1}+...+h_j$ is central in $\hgl_n$ if $n|j-i$). In order to prove \eqref{eqn:cohom} by induction on $i-j$, it suffices to replace all $(\of_{[a+1;j+1)} + \of_{[a+1;j)})\tA_m$ in the LHS of \eqref{eqn:two line} by $\tA_m(\of_{[a+1;j)} + \of_{[a;j)} + m \alpha_{j-a-1} )$ and show that \eqref{eqn:two line} is an identity. Upon making this replacement, the left-hand side of \eqref{eqn:two line} equals: 
\begin{multline*}
\tA_m \left[ \sum_{a = i}^{j-1} (-1)^{i-a} \orr_{\left \lfloor \frac {a-i}n \right \rfloor} \left(\of_{[a+1;j)} + \of_{[a;j)} + m \alpha_{j-a-1} \right) \right] = \\ = \tA_m \left[ \of_{[i;j)} + \sum_{a=i+1}^{j-1} (-1)^{i-a} \left(\orr_{\left \lfloor \frac {a-i}n \right \rfloor} - \orr_{\left \lfloor \frac {a-i-1}n \right \rfloor} \right) \of_{[a;j)}  + (-1)^{i-j+1} \orr_{\left \lfloor \frac {j-1-i}n \right \rfloor} h_j + \right. \\ \left. + m \sum_{a=i}^{j-1} (-1)^{i-a} \orr_{\left \lfloor \frac {a-i}n \right \rfloor} \alpha_{j-a-1} \right] =  \tA_m \left[ \of_{[i;j)} + \sum^{n|a-i}_{i < a < j}  \frac {(-1)^{i-a} \orr_{\left \lfloor \frac {a-i-1}n \right \rfloor}\of_{[a;j)} }{p\left \lfloor \frac {a-i}n \right \rfloor} + \right. \\ \left.  + (-1)^{i-j+1} \orr_{\left \lfloor \frac {j-1-i}n \right \rfloor} h_j + m \sum_{a=i}^{j-1} (-1)^{i-a} \orr_{\left \lfloor \frac {a-i}n \right \rfloor} \alpha_{j-a-1} \right]
\end{multline*}
The formula above matches the right-hand side of \eqref{eqn:two line}, as we needed to show, as a consequence of the combinatorial identity:
$$
\sum_{a=i}^{j-1} (-1)^{i-a} \orr_{\left \lfloor \frac {a-i}n \right \rfloor} \alpha_{j-a-1} = (-1)^{i-j+1} \left[ \orr_{\left \lfloor \frac {j-i}n \right \rfloor} (j-i) - \orr_{\left \lfloor \frac {j-i-1}n \right \rfloor} (j-i-1) \right]
$$ 
(the identity above is a simple exercise that follows from $\orr_k = {k + \frac 1p \choose k }$, which we leave to the interested reader). 

\end{proof}

\subsection{} Because $H$ is isomorphic to the Verma module $M$ via \eqref{eqn:iso cohom}, we may regard $\tA_m$ as an operator $M' \rightarrow M$ between Verma modules with (possibly different) universal highest weights $a_1',...,a_n'$ and $a_1,...,a_n$, respectively. It is easy to see that Proposition \ref{prop:unique} has a cohomological analogue, meaning that there is at most a unique (up to scalar multiple) operator satisfying the commutation relations \eqref{eqn:cohom}. \\

\begin{proposition}
\label{prop:m=0}

If $M = M'$ (i.e. $a_i=a_i'\ \forall i$), the operator $A_0 \Big|_{M = M'}$ acts as:
\begin{equation}
\label{eqn:big}
g := \begin{pmatrix} 1 & 0 & 0 & z^{-1} \\ 1 & 1 & 0 & 0 \\ 0 & 1 & 1 & 0 \\ 0 & 0 & 1 & 1 \end{pmatrix}
\end{equation}
which is an element $\in \widehat{GL}_n^- \hookrightarrow \widehat{U(\hgl_n^-)}$ explicitly given by the formula:
\begin{equation}
\label{eqn:big formula}
g = \exp \left( \sum_{1\leq i \leq n}^{j>i} \frac {(-1)^{j-i-1}}{j-i} \cdot E_{ji} \right)
\end{equation}
Any power series in $\{E_{ji}\}_{i<j}$ acts correctly on the Verma module $M$. \\
\end{proposition}

\begin{proof} Let us first prove that the right-hand sides of \eqref{eqn:big} and \eqref{eqn:big formula} are equal. Letting $S = g - \text{Id} = \sum_{i=1}^n E_{i+1,i}$, then we observe that:
$$
S^k = \sum_{i=1}^n E_{i+k,i}
$$
$\forall k \geq 1$, and the required identity is a consequence of the Taylor series expansion: 
$$
\text{Id} + S = \exp\left(\sum_{k=1}^\infty \frac {(-1)^{k-1}}k \cdot S^k \right)
$$
Let us now show that $A_0|_{M = M'} = g$. Because of the uniqueness statement immediately preceding the statement of the Proposition, it suffices to show that:
\begin{equation}
\label{eqn:want equality}
(E_{ij} + E_{i+1,j}) g = g (E_{ij}+ E_{i,j-1})
\end{equation}
$\forall i,j \in \BZ$. The equality \eqref{eqn:want equality} may be regarded as the adjoint action of $g \in GL_n[z^{-1}]$ on elements of the Lie algebra $\fgl_n[z^{\pm 1}]$, and therefore it suffices to check that the equality holds in the space of $n \times n$ matrices with coefficients Laurent polynomials in $z$. In this case, formula \eqref{eqn:want equality} reduces to the following easily checked matrix identity (seen below for $n=4$, $i=3$, $j=2$):
$$
\begin{pmatrix} 0 & 0 & 0 & 0 \\ 0 & 0 & 0 & 0 \\ 0 & 1 & 0 & 0 \\ 0 & 1 & 0 & 0 \end{pmatrix} \begin{pmatrix} 1 & 0 & 0 & z^{-1} \\ 1 & 1 & 0 & 0 \\ 0 & 1 & 1 & 0 \\ 0 & 0 & 1 & 1 \end{pmatrix} = \begin{pmatrix} 1 & 0 & 0 & z^{-1} \\ 1 & 1 & 0 & 0 \\ 0 & 1 & 1 & 0 \\ 0 & 0 & 1 & 1 \end{pmatrix} \begin{pmatrix} 0 & 0 & 0 & 0 \\ 0 & 0 & 0 & 0 \\ 1 & 1 & 0 & 0 \\ 0 & 0 & 0 & 0 \end{pmatrix}
$$

\end{proof}

\subsection{} Let $V$ be a vector space with basis given by symbols $y_1^{b_1}...y_n^{b_n}$, as the exponents $b_1,...,b_n$ range over all complex numbers. There exists an action:
\begin{equation}
\label{eqn:v infinite}
\hgl_n \curvearrowright V
\end{equation}
where the central charge $\gamma$ acts by 0, and:
\begin{equation}
\label{eqn:finite action}
E_{ij} \leadsto (-1)^{i-j} y_{\bari} \frac {\partial}{\partial y_{\barj}} - \text{Id} \cdot  m \begin{cases} 1 & \text{if } i=j \\ (-1)^{i-j-1}\frac np & \text{if } i - j \in n \BN \\ 0 & \text{otherwise }  \end{cases} 
\end{equation}
(recall that we write $\bari$ for the residue class of $i$ in the set $\{1,...,n\}$). \\


\begin{theorem}
\label{thm:rep}

The operator $M' \stackrel{\tA_m}\longrightarrow M$ factors as:
\begin{equation}
\label{eqn:comp 0}
\tA_m : M' \stackrel{g}\longrightarrow M' \stackrel{\Phi_m}\longrightarrow M \otimes V \stackrel{\emph{Id} \otimes \emph{ev}}\longrightarrow M
\end{equation}
where $M' \stackrel{\Phi_m}\longrightarrow M \otimes V$ is the $\hgl_n$ intertwiner, and $V \stackrel{\emph{ev}}\longrightarrow \BC$ denotes the linear map given by evaluating polynomials at $y_1 = ... = y_n = 1$. \\

 		
\end{theorem}

\noindent We refer to $\Phi_m$ as ``the" intertwiner because it is unique up to constant multiple, by a straightforward analogue of Proposition \ref{prop:unique}. This crucially uses the fact that, in order for it to commute with $E_{11},...,E_{nn}$, the operator $\Phi_m$ must map the generators of the Verma module according to:
\begin{equation}
\label{eqn:inter coeffs}
v' \stackrel{\Phi_m}\leadsto \alpha v \otimes \left( y_1^{a_1'-a_1+m}... y_n^{a_n'-a_n+m} \right) + \text{higher order}
\end{equation}
for some scalar $\alpha$. Just like at the end of the proof of Proposition \ref{prop:unique}, knowledge of the constant $\alpha$ allows one to iteratively compute all higher order terms in the expression above, which ultimately determines $\Phi_m$ completely. \\


\begin{proof} If $\BC$ is given the trivial $\hgl_n$ module structure, we have $(E_{ij} + E_{i+1,j}) \circ \ev = 0$ for all $i,j \in \BZ$. Meanwhile, because of formulas \eqref{eqn:finite action}, it is easy to see that:
$$
\ev \circ (E_{ij} + E_{i+1,j}) = - m \alpha_{i-j} \cdot \ev
$$
where $\alpha_{i-j}$ is defined in \eqref{eqn:alpha}. By taking the difference, we conclude that:
\begin{equation}
\label{eqn:hari}
[E_{ij}+E_{i+1,j}, \ev] = m \alpha_{i-j} \cdot \ev 
\end{equation}
Let us now prove that $\tA_m$ factors as stipulated in \eqref{eqn:comp}. Because of (the Lie algebra analogue of) Proposition \ref{prop:unique}, it is enough to show that the composition of operators in \eqref{eqn:comp 0} satisfies the commutation relations \eqref{eqn:cohom} in place of $\tA_m$. Explicitly, this is because:
\begin{multline*}
(E_{ij} + E_{i+1,j}) \circ (\text{Id} \otimes \ev) \Phi_m g \stackrel{\eqref{eqn:hari}}= (\text{Id} \otimes \ev) \circ (E_{ij} + E_{i+1,j} + m \alpha_{i-j}) \Phi_m g = \\
= (\text{Id} \otimes \ev) \Phi_m  \circ (E_{ij} + E_{i+1,j} + m \alpha_{i-j}) g \stackrel{\eqref{eqn:want equality}}= (\text{Id} \otimes \ev) \Phi_m g (E_{ij} + E_{i,j-1} + m \alpha_{i-j})
\end{multline*}
where the middle equality holds on account of $\Phi_m$ being an intertwiner. 

\end{proof} 

\subsection{}

For the remainder of the present Section, we assume $m \in \BN$ and $M = M'$, by which we mean that we equate the equivariant parameters $a_i = a_i'$. In this case, the representation $V$ of \eqref{eqn:v infinite} can be replaced by the finite-dimensional representation:
\begin{equation}
\label{eqn:symm}
\hgl_n \curvearrowright S^{mn} = \BC[y_1,...,y_n]^{\text{total degree }mn}
\end{equation}
with the action given by formula \eqref{eqn:finite action}. It is easy to observe that $S^{mn}$ is simply the $mn$--th symmetric power of $\BC^n$ endowed with the following action (fix a basis $y_1,...,y_n$ of $\BC^n$, although strictly speaking it should be the dual basis to \eqref{eqn:symm}):
\begin{equation}
\label{eqn:taut action}
E_{ij} \cdot y_k = (-1)^{i-j} \delta_k^{\barj} y_{\bari} - y_k \begin{cases} \frac 1n & \text{if } i = j \\ (-1)^{i-j-1} \frac 1p & \text{if } i - j \in n \BN \\ 0 & \text{otherwise }  \end{cases} 
\end{equation}
$\forall i,j\in \BZ$. In other words, $\BC^n$ is the evaluation representation corresponding to the tautological representation of $\fgl_n$ and $z = 1$, twisted by a character of the central Heisenberg subalgebra. Because $a_i = a_i'$, the lowest weight coefficient $(y_1...y_n)^m$ from \eqref{eqn:inter coeffs} actually lies in $S^{mn}$, which implies that:
\begin{equation}
\label{eqn:comp}
\tA_m : M \stackrel{g}\longrightarrow M \stackrel{\Phi_m}\longrightarrow M \otimes S^{mn} \stackrel{\text{Id} \otimes \ev}\longrightarrow M
\end{equation}
where $\Phi_m$ is the $\hgl_n$ intertwiner (the proof of \eqref{eqn:comp} is identical to that of \eqref{eqn:comp 0}). \\

\subsection{} 

Recall that the Verma module has a grading:
$$
M = \bigoplus_{\bd = (d_1,...,d_n) \in \nn} M_{\bd}
$$
Let us consider the diagonal operator:
\begin{equation}
\label{eqn:grading}
\bx^{\bd} : M \longrightarrow M
\end{equation}
which acts on each $M_{\bd}$ as multiplication by $x_1^{d_1}...x_n^{d_n}$, where $x_1,...,x_n$ are formal parameters. The generalized character of $\Phi_m : M \rightarrow M \otimes S^{mn}$ is defined as:
$$
\chi_{\Phi_m} = \text{Tr} \left(\Phi_m \bx^{\bd} \right) \in S^{mn}[x_1,...,x_n]
$$
Note that on virtue of $\Phi_m$ being an intertwiner, the character $\chi_{\Phi_m}$ actually takes values in $(y_1...y_n)^m \cdot \BC[[x_1,...,x_n]]$, because the weight zero subspace of $S^{mn}$ is the span of the vector $(y_1...y_n)^m$. The following Proposition closely follows \cite{EK}. \\

\begin{proposition}
\label{prop:char inter}

The character $\chi_{\Phi_m}$ satisfies the identity:
\begin{equation}
\label{eqn:char inter}
\CO_m \cdot \chi_{\Phi_m} = 0
\end{equation}
where $\CO_m$ is the following differential operator on $\BC[[x_1,...,x_n]]$:
\begin{equation}
\label{eqn:mathcal}
\CO_m = \sum_{i=1}^n D_i(D_i - D_{i+1} + a_i - a_{i+1}) - \sum_{k=1}^\infty \frac {kn(m+1)}{1-x^{-k}} + 
\end{equation}
$$
+ \sum_{(i,j) \in \frac {\BZ^2}{(n,n)\BZ}}^{i<j} \left[ \frac {a_i - a_j - i + j + D_i - D_{i-1} - D_j + D_{j-1}}{1-x_{[i;j)}^{-1}} + \frac {m(m+1)}{(1-x_{[i;j)})(1-x_{[i;j)}^{-1})} \right] 
$$
where $x = x_1...x_n$ and $x_{[i;j)} = x_i...x_{j-1}$. \\
	
\end{proposition}

\begin{proof} Consider the following element in the completion $\widehat{U(\hgl_n)}$:
$$
c =  \sum_{i=1}^n \left( \frac {E_{ii}^2}2 + \frac {\gamma i}n E_{ii} \right)  + \sum_{(i,j) \in \frac {\BZ^2}{(n,n)\BZ}}^{i<j} E_{ij} E_{ji} - \frac 1{\gamma} \sum_{k=1}^\infty P_kP_{-k}
$$ 
Let us note that:
\begin{equation}
\label{eqn:comm}
[c,E_{ab}] = \left( 1-\frac {\gamma}n \right) (b-a) E_{ab}
\end{equation}
as follows by adding the identities below (assume $a < b$ for conciseness):
$$
\sum^{i<j}_{(i,j) \in \frac {\BZ^2}{(n,n)\BZ}} \left[ E_{ij} E_{ji},E_{ab} \right] =  \sum_{(i,j) \in \frac {\BZ^2}{(n,n)\BZ}}^{i<j} \Big( E_{ij} [E_{ji}, E_{ab}] + [E_{ij}, E_{ab}] E_{ji} \Big) = 
$$
$$
= \sum_{j > a} E_{aj}E_{jb} - \sum_{i < b} E_{ib}E_{ai} + \sum_{i < a} E_{ib} E_{ai} - \sum_{j > b} E_{aj} E_{jb} =  \sum_{j=a+1}^{b} E_{aj} E_{jb} - \sum_{i=a}^{b-1} E_{ib}E_{ai} = 
$$
$$
= E_{bb} E_{ab} - E_{ab} E_{aa} + \sum_{i=a+1}^{b} [E_{ai}, E_{ib}] =  E_{bb} E_{ab} - E_{ab} E_{aa} + (b - a) E_{ab} - \delta \cdot \frac {b-a}n P_{\frac {b-a}n}
$$
(the symbol $\delta$ denotes 1 if $n|a-b$ and 0 otherwise) together with:
$$
\sum_{i=1}^n \left[ \frac {E_{ii}^2}2 + \frac {\gamma i}n E_{ii}, E_{ab} \right] = \sum_{i=1}^n \left( \frac {E_{ii}}2 [E_{ii}, E_{ab}] + [E_{ii}, E_{ab}] \frac{E_{ii}}2 + \frac {\gamma i}n [E_{ii}, E_{ab}] \right) = 
$$
$$
= E_{ab} E_{aa} - E_{bb} E_{ab} + \gamma \left( \left \lfloor \frac {a-1}n \right \rfloor - \left \lfloor \frac {b-1}n \right \rfloor \right) E_{ab} + \gamma \frac {\bar{a}-\bar{b}}n E_{ab}
$$
(above we used the relation $E_{aa} = E_{\bar{a}\bar{a}} - \gamma \left \lfloor \frac {a-1}n \right \rfloor$, where $\bar{a}$ denotes the residue class of the integer $a$ in the set $\{1,...,n\}$) and:
$$
- \frac 1\gamma \sum_{k=1}^\infty [P_k P_{-k}, E_{ab}] = \delta \cdot \frac {b-a}n P_{\frac {b-a}n}
$$
If we consider the degree operator:
$$
\deg :  M \rightarrow M
$$
which acts on $M_{\bd}$ as multiplication by $d_1+...+d_n$, then formula \eqref{eqn:comm} implies:
$$
c \Big|_M = - \frac pn \deg + c \Big|_{M_0} = -\frac pn \deg + \sum_{i=1}^n \left[ \frac {(a_i-i)^2}2 + \frac {(n+p) i (a_i-i)}n \right]
$$
Since $\chi_{\Phi_m \deg} = D \chi_{\Phi_m}$, where $D = D_1+...+D_n$, we obtain the identity:
\begin{equation}
\label{eqn:identity 1}
\chi_{\Phi_m c} = \left( -\frac pn D + \sum_{i=1}^n \left[ \frac {(a_i-i)^2}2 + \frac {(n+p) i (a_i-i)}n \right] \right) \chi_{\Phi_m}
\end{equation}
On the other hand, because $E_{ii}$ acts on $M_{\bd}$ as multiplication by $a_i - i + d_i - d_{i-1}$:
\begin{equation}
\label{eqn:identity 2}
\chi_{\Phi_m E_{ii}} = \Tr(\Phi_m E_{ii} \bx^{\bd}) = (a_i - i + D_i - D_{i-1}) \chi_{\Phi_m}
\end{equation}
and analogously:
\begin{equation}
\label{eqn:identity 3}
\chi_{\Phi_m E_{ii}^2} = \Tr(\Phi_m E_{ii}^2 \bx^{\bd}) = (a_i - i + D_i - D_{i-1})^2 \chi_{\Phi_m}
\end{equation}
Moreover, if $i<j$ we have:
$$
\Tr(\Phi_m E_{ij} \bx^{\bd}) = \frac {\Tr(\Phi_m \bx^{\bd} E_{ij})}{x_{[i;j)}} = \frac {(E_{ij} \otimes 1) \Tr(\Phi_m \bx^{\bd})}{x_{[i;j)}} =
$$
\begin{equation}
\label{eqn:identity 4}
= \frac {- E_{ij} \Tr(\Phi_m \bx^{\bd})}{x_{[i;j)}} + \frac {\Tr(\Phi_m E_{ij} \bx^{\bd})}{x_{[i;j)}} \ \Rightarrow \ \chi_{\Phi_m E_{ij}} = \frac {E_{ij} \cdot \chi_{\Phi_m}}{1-x_{[i;j)}}
\end{equation}
where the third equality follows from $(E_{ij} \otimes 1 + 1 \otimes E_{ij}) \Phi_m = \Phi_m E_{ij}$, namely the intertwiner property of the operator $\Phi_m$. Analogously, we have:
\begin{equation}
\label{eqn:identity 5}
\chi_{\Phi_m P_k} = \frac {P_k \cdot \chi_{\Phi_m}}{1-x^k}
\end{equation}
(recall that $x_{[i;j)} = x_i..x_{j-1}$ and $x = x_1...x_n$). Moreover, we have:
$$
\Tr(\Phi_m E_{ij} E_{ji} \bx^{\bd}) = x_{[i;j)} \Tr(\Phi_m E_{ij} \bx^{\bd} E_{ji}) = x_{[i;j)} \Tr((E_{ji} \otimes 1) \Phi_m E_{ij} \bx^{\bd}) =
$$
$$
= -x_{[i;j)} E_{ji} \cdot \Tr(\Phi_m E_{ij} \bx^{\bd}) + x_{[i;j)} \Tr(\Phi_m E_{ji} E_{ij} \bx^{\bd}) = - x_{[i;j)} \frac {E_{ji} E_{ij} \cdot \chi_{\Phi_m}}{1-x_{[i;j)}} +
$$
$$
+ x_{[i;j)} \Tr(\Phi_m E_{ij} E_{ji} \bx^{\bd}) + x_{[i;j)} \Tr(\Phi_m (E_{jj} - E_{ii}) \bx^{\bd}) \quad \Rightarrow \quad \chi_{\Phi_m E_{ij}E_{ji}} =
$$
\begin{equation}
\label{eqn:identity 6}
= \frac {E_{ji} E_{ij} \chi_{\Phi_m}}{(1-x_{[i;j)})(1-x_{[i;j)}^{-1})} + \frac {(a_i - a_j - i + j + D_i - D_{i-1} - D_j + D_{j-1})\chi_{\Phi_m}}{1-x_{[i;j)}^{-1}}
\end{equation}
Note that we may replace $E_{ji}E_{ij}$ by $E_{ij}E_{ji}$ on the last row, since the generalized character $\chi_{\Phi_m}$ having degree 0 implies that $E_{ii} \chi_{\Phi_m} = 0$. Analogously, we have:
\begin{equation}
\label{eqn:identity 7}
\chi_{\Phi_m P_k P_{-k}} =  \frac {P_kP_{-k} \chi_{\Phi_m}}{(1-x^k)(1-x^{-k})} + \frac {kn(n+p) \chi_{\Phi_m}}{1-x^{-k}}
\end{equation}
Taking appropriate linear combinations of equations \eqref{eqn:identity 1}, \eqref{eqn:identity 2}, \eqref{eqn:identity 3}, \eqref{eqn:identity 6} and \eqref{eqn:identity 7}, we obtain the following identity:
$$
\chi_{\Phi_m c} =  \sum_{i=1}^n \left[ \frac {\text{RHS of \eqref{eqn:identity 3}}}2 + i \left(1 + \frac pn \right) \Big( \text{RHS of \eqref{eqn:identity 2}} \Big) \right] + 
$$
\begin{equation}
\label{eqn:identity 8}
+ \sum_{(i,j) \in \frac {\BZ^2}{(n,n)\BZ}}^{i<j} \Big( \text{RHS of \eqref{eqn:identity 6}} \Big) - \frac 1{n+p} \sum_{k=1}^\infty \Big( \text{RHS of \eqref{eqn:identity 7}} \Big)
\end{equation}
If we recall the action \eqref{eqn:finite action} and the fact that $\chi_{\Phi_m}$ equals $(y_1...y_n)^m$ times a power series in the variables $x_1,...,x_n$, then we observe that:
\begin{align*}
&E_{ij}E_{ji} \cdot \chi_{\Phi_m} = m(m+1) \chi_{\Phi_m} \qquad \quad \text{if } \bari \neq \barj \\
&E_{ij}E_{ji} \cdot \chi_{\Phi_m} = m^2\left(1+\frac np \right) \chi_{\Phi_m} \ \ \quad \text{if } \bari = \barj \\
&P_kP_{-k} \cdot \chi_{\Phi_m} = mn\left(mn+\frac {n^2}p \right) \chi_{\Phi_m}
\end{align*}
With this in mind, as well as the straightforward identity:
\begin{equation}
\label{eqn:mich}
\sum_{k=1}^{\infty} \frac k{1-x^{-k}} = \sum_{k=1}^{\infty} \frac 1{(1-x^k)(1-x^{-k})} \quad \left( = - \sum_{k,l=1}^\infty k x^{kl} \right)
\end{equation}
setting \eqref{eqn:identity 1} equal to \eqref{eqn:identity 8} yields precisely \eqref{eqn:char inter}.

\end{proof}

\subsection{} 

It is easy to observe that $\hgl_n$ acts faithfully on the representation:
\begin{equation}
\label{eqn:infinite rep}
\BC^\infty = \BC^n[z,z^{-1}]
\end{equation}
with $\gamma \mapsto 0$. We will fix a basis $\{y_k\}_{k \in \BZ}$ of \eqref{eqn:infinite rep}, with the understanding that:
$$
y_k = y_{\bark} \cdot z^{-  \left \lceil \frac {k-1}n \right \rceil}
$$
In terms of this basis and the notation \eqref{eqn:notation}, the action $\hgl_n \curvearrowright \BC^\infty$ is given by:
$$
E_{ij} \cdot y_k = \delta_{\barj}^{\bark} y_{k-j+i}, \qquad \forall i,j,k \in \BZ
$$
It is easy to see that we have the following multiplication rules in $\BC^\infty$:
\begin{equation}
\label{eqn:rules}
E_{ij} \cdot E_{i'j'} = \delta_{\barj}^{\barii} E_{i,j+j'-i'}
\end{equation}
Note that elements in $\widehat{GL}_n^-$ are well-defined in the representation \eqref{eqn:infinite rep}, if one is willing to allow infinite sums of $y_k$'s in the direction of positive $k$. Consider:
$$
\bx^{\bd} : \BC^\infty \longrightarrow \BC^\infty, \qquad \bx^{\bd} \cdot y_k = y_k \begin{cases} x_1...x_{k-1} &\text{if } k>0 \\ \frac 1{x_k...x_0} &\text{if } k \leq 0 \end{cases}
$$
where $x_1,...,x_n$ are the same formal symbols as before. It is no coincidence that we use the same symbol for the operator above as for the endomorphism \eqref{eqn:grading}, since they both have the same commutation rules with the affine Lie algebra generators:
\begin{align}
&\bx^{\bd} \cdot E_{ij} \cdot \bx^{-\bd} = x_i...x_{j-1} \cdot E_{ij} \label{eqn:om 1} \\
&\bx^{\bd} \cdot E_{ji} \cdot \bx^{-\bd} = \frac 1{x_i...x_{j-1}} \cdot E_{ji} \label{eqn:om 2}
\end{align}
for any $i < j \in \BZ$. \\

\begin{proposition}
\label{prop:identity}

We have the following identity in the representation $\BC^\infty$:
\begin{equation}
\label{eqn:identity}
g \cdot (\emph{\bx}^{\bd} \cdot h \cdot {\emph{\bx}}^{-\bd}) = h
\end{equation}
where $h \in \widehat{GL}_n^-[[x_1,...,x_n]]$ acts in the representation $\BC^\infty$ as:
\begin{equation}
\label{eqn:formula h}
h = 1 + \sum^{i < j}_{(i,j) \in \BZ^2/(n,n)\BZ} E_{ji} \prod_{a=i}^{j-1} \frac {x_a}{x_i...x_a - 1}
\end{equation}

\end{proposition}

\begin{proof} By \eqref{eqn:om 2}, we have:
$$
\bx^{\bd} \cdot h \cdot \bx^{-\bd} = 1 + \sum^{i < j}_{(i,j) \in \BZ^2/(n,n)\BZ} E_{ji} \prod_{a=i}^{j-1} \frac 1{x_i...x_a - 1}
$$
Since $g = 1 + \sum_{i=1}^n E_{i+1,i}$ in the representation $\BC^\infty$ according to \eqref{eqn:big}, identity \eqref{eqn:identity} is an easy consequence of \eqref{eqn:rules}. 

\end{proof}

\subsection{} Since $g$, $h$ and $\bx^{\bd} h \bx^{-\bd}$ are all elements of $\widehat{GL}_n^-[[x_1,...,x_n]]$, the fact that relation \eqref{eqn:identity} holds in the faithful representation $\BC^\infty$ implies that it also holds in the Verma module $M$. Recall the Weyl denominator $\delta_n$ of \eqref{eqn:weyl}, and consider also:
\begin{equation}
\label{eqn:weyl one}
\delta_1 = \prod_{k=1}^\infty (1-x^k)
\end{equation}
where $x = x_1...x_n$. \\

\begin{proposition}
\label{prop:characters}

If $M \stackrel{\Phi_m}\longrightarrow M \otimes S^{mn}$ is the $\hgl_n$--intertwiner, then we have:
\begin{equation}
\label{eqn:ev}
\emph{\ev} \circ \chi_{\Phi_m g} = \chi_{\Phi_m} \Big|_{y_1,...,y_n \mapsto 1} \cdot \delta_n^{-m} \delta_1^{-\frac {mn}p}
\end{equation}
(recall that $\chi_{\Phi_m} \in (y_1...y_n)^m \cdot \BC[[x_1,...,x_n]]$). The LHS equals $\chi_{\tA_m}$ by \eqref{eqn:comp}. \\

\end{proposition}

\begin{proof} Note that:
\begin{multline}
\label{eqn:displayed}
\chi_{\Phi_m g} = \Tr(\Phi_m g \bx^{\bd}) \stackrel{\eqref{eqn:identity}}= \Tr(\Phi_m h \bx^{\bd} h^{-1}) = \\
= \Tr((h \otimes h) \Phi_m \bx^{\bd} h^{-1}) = h \cdot \Tr(\Phi_m \bx^{\bd}) = h \cdot \chi_{\Phi_m}
\end{multline}
Therefore, we must compute how $h \in \widehat{GL}_n^-[[x_1,...,x_n]]$ acts on $(y_1...y_n)^m \in S^{mn}$. To this end, we note that the map:
\begin{equation}
\label{eqn:assignment}
\BC^\infty \rightarrow \BC^n \otimes \BC, \qquad y_i \mapsto y_{\bari} \otimes 1 
\end{equation}
is a $\hgl_n$ intertwiner if $\BC^n$ denotes the representation \eqref{eqn:taut action} and $\BC$ is the representation corresponding to the character:
$$
A \cdot z^k \leadsto \begin{cases} 0 &\text{if } k>0 \\ \frac 1n &\text{if } k = 0 \\ \frac{(-1)^{nk-1} \text{Tr}(A)}p &\text{if }k<0 \end{cases}
$$
Since $h$ is a group-like element, the assignment \eqref{eqn:assignment} sends:
\begin{equation}
\label{eqn:assignment 2}
(h \cdot y_i) \leadsto (h \cdot y_{\bari}) \otimes (h \cdot 1)
\end{equation} 
The left-hand side can be computed by using the explicit presentation \eqref{eqn:formula h}:
\begin{equation}
\label{eqn:assignment 3}
h \cdot y_i = y_i + \sum_{i < j} (-1)^{i-j} y_j \prod_{a=i}^{j-1} \frac {x_a}{x_i...x_a - 1}
\end{equation}
Meanwhile, using the decomposition $\hgl_n = \hsl_n \oplus \hgl_1$, we may write $g$ and $h$ as:
\begin{align*}
&g = g_{\hsl_n} \cdot g_{\hgl_1} \\
&h = h_{\hsl_n} \cdot h_{\hgl_1}
\end{align*}
and the identity \eqref{eqn:identity} holds component-wise in terms of the factorization above. Therefore, formula \eqref{eqn:big formula} implies that:
$$
g_{\hgl_1} = \exp \left( \sum_{k=1}^\infty \frac {(-1)^{nk-1} P_{-k}}{nk} \right) \quad \stackrel{\eqref{eqn:identity}}\Longrightarrow \quad h_{\hgl_1} = \exp \left( \sum_{k=1}^\infty \frac {(-1)^{nk-1} P_{-k}}{nk(1-x^{-k})} \right) 
$$
In \eqref{eqn:assignment 2}, the action of $h$ on the generator $1 \in \BC$ factors through $h_{\hgl_1}$, hence:
\begin{equation}
\label{eqn:assignment 4}
h \cdot 1 = h_{\hgl_1} \cdot 1  = \exp \left(\sum_{k=1}^\infty \frac {1}{pk(1-x^{-k})} \right) = \prod_{k=1}^\infty (1-x^k)^{\frac 1p}
\end{equation}
If we plug \eqref{eqn:assignment 3} and \eqref{eqn:assignment 4} into \eqref{eqn:assignment 2} and apply the linear map $\ev$ which sends all $y_i$'s to 1, then we conclude that:
$$
1 + \sum_{i < j} \prod_{a=i}^{j-1} \frac {(-1)^{j-i}  x_a}{x_i...x_a - 1} = \ev(h\cdot y_{\bari}) \cdot \prod_{k=1}^\infty (1-x^k)^{\frac 1p} \quad \Rightarrow
$$
$$
\Rightarrow \quad \ev(h\cdot y_{\bari}) = \prod_{i<j}  \frac 1{1-x_i...x_{j-1}} \cdot \prod_{k=1}^\infty (1-x^k)^{-\frac 1p}
$$
where we use the identity $\prod_{a=i}^{j-1} \frac {(-1)^{j-i} x_a}{x_i...x_a - 1} = \prod_{a=i}^{j-1} \frac 1{1 - x_i...x_a} - \prod_{a=i}^{j-2} \frac 1{1 - x_i...x_a}$.  Since $S^{mn}$ is the $mn$-th symmetric power of the representation $\BC^n$ and $h$ is group-like:
$$
\ev(h\cdot (y_1...y_n)^m) = \left( \prod^{i < j}_{(i,j) \in \BZ^2/(n,n)\BZ} \frac 1{1-x_i...x_{j-1}} \right)^m \prod_{k=1}^\infty (1-x^k)^{-\frac {mn}p} = \delta_n^{-m} \cdot \delta_1^{-\frac {mn}p}
$$
According to \eqref{eqn:displayed}, this implies \eqref{eqn:ev}. 

\end{proof}

\begin{proof} \textbf{of Theorem \ref{thm:main}:} Let us define $Y_m$ as in \eqref{eqn:finally equal}, and the task becomes to prove that it satisfies the differential equation \eqref{eqn:function}. Since the coefficients of this differential equation (in the Taylor series expansion with respect to $x_1,...,x_n$) are all polynomials in $m$, it suffices to prove that \eqref{eqn:function} holds in the case $m \in \BN$. Therefore, the cohomological version of Proposition \ref{prop:char}, together with \eqref{eqn:ev}, imply that:
$$
Y_m = \chi_{\Phi_m} \Big|_{y_1,...,y_n \mapsto 1} \cdot x_1^{b_1}...x_n^{b_n} \delta_n \delta_1^{-\frac {mn}p}
$$
hence \eqref{eqn:char inter} implies:
$$
\CO_m \left(Y_m \cdot x_1^{-b_1}...x_n^{-b_n} \delta^{-1}_n \delta_1^{\frac {mn}p} \right) = 0
$$
where $\CO_m$ is the differential operator \eqref{eqn:mathcal}. By the Leibniz rule, we obtain:
$$
0 = \CO_m Y_m + \sum_{i=1}^n (D_i-D_{i-1}) Y_m \cdot (D_i-D_{i-1}) \left( x_1^{-b_1}...x_n^{-b_n} \delta^{-1}_n \delta_1^{\frac {mn}p} \right) + 
$$
$$
+ Y_m \sum_{i=1}^n \frac {(D_i - D_{i-1})^2}2 \left( x_1^{-b_1}...x_n^{-b_n} \delta^{-1}_n \delta_1^{\frac {mn}p} \right) +
$$
$$
+ Y_m \left[ \left(a_i + \frac {ip}n \right) (D_i - D_{i-1}) + \frac pn D \right]\left( x_1^{-b_1}...x_n^{-b_n} \delta^{-1}_n \delta_1^{\frac {mn}p} \right)
$$
\begin{equation}
\label{eqn:leibniz}
+ Y_m \sum_{(i,j) \in \frac {\BZ^2}{(n,n)\BZ}}^{i<j} \frac {(D_i-D_{i-1}-D_j+D_{j-1})  \left( x_1^{-b_1}...x_n^{-b_n} \delta^{-1}_n \delta_1^{\frac {mn}p} \right)}{1-x_{[i;j)}^{-1}} \qquad
\end{equation}
It is straightforward to obtain \eqref{eqn:function} from the equality above, by using the identities:
\begin{align}
&(D_s - D_{s-1}) \cdot (x_1^{-b_1}...x_n^{-b_n}) = (b_{s-1} - b_s) (x_1^{-b_1}...x_n^{-b_n}) \label{eqn:has 1} \\
&(D_s - D_{s-1}) \cdot \delta_1^{\frac {mn}p} = 0 \label{eqn:has 2} \\
&D \cdot (x_1^{-b_1}...x_n^{-b_n}) = -(b_1+...+b_n) (x_1^{-b_1}...x_n^{-b_n}) \label{eqn:has 3} \\
&D \cdot \delta_1^{\frac {mn}p} = \frac {mn}p \delta_1^{\frac {mn}p} \sum_{k=1}^\infty \frac {nk}{1-x^{-k}} \label{eqn:has 4}
\end{align}
as well as:
$$
D_s \delta^{-1}_n = \delta^{-1}_n \sum^{i < j}_{(i,j) \in \BZ^2/(n,n)\BZ} \frac {x_i...x_{j-1}}{1-x_i...x_{j-1}} \Big( \# \{i,...,j-1\} \cap \{s + n \BZ\} \Big)
$$
which, in turn, implies the formulas:
\begin{align}
&D \delta_n^{-1} = \delta_n^{-1} \sum^{i < j}_{(i,j) \in \BZ^2/(n,n)\BZ} \frac {i-j}{1-x_{[i;j)}^{-1}} \label{eqn:has 5} \\
&(D_s - D_{s-1})\delta^{-1}_n = \delta^{-1}_n \rho_s \label{eqn:has 6}
\end{align}
$\forall s$, where:
$$
\rho_s =  \sum_{i < s} \frac 1{1-x_{[i;s)}^{-1}} - \sum_{s<j} \frac 1{1-x_{[s;j)}^{-1}} 
$$
Applying $D_s$ to \eqref{eqn:has 6} yields:
$$
D_s(D_s-D_{s-1})\delta^{-1}_n = \delta^{-1}_n \left[\sum_{i<s} \frac {\left \lfloor \frac {s-i}n \right \rfloor}{(1-x_{[i;s)})(1-x_{[i;s)}^{-1})} - \sum_{s<j} \frac {\left \lceil \frac {j-s}n \right \rceil }{(1-x_{[s;j)})(1-x_{[s;j)}^{-1})} \right]
$$
$$
- \delta_n^{-1} \sum_{s=1}^n \sum^{i < j}_{(i,j) \in \BZ^2/(n,n)\BZ} \frac {\# \{i,...,j-1\} \cap \{s + n \BZ\} }{1-x_{[i;j)}^{-1}} \cdot \rho_s
$$
Summing the equation above over all $s \in \{1,...,n\}$ implies:
$$
\sum_{s=1}^n \frac {(D_s-D_{s-1})^2}2 \delta^{-1}_n = \delta_n^{-1}  \sum^{i < j}_{(i,j) \in \BZ^2/(n,n)\BZ} \left[ \frac {\delta_{\bari}^{\barj}-1}{(1-x_{[i;j)})(1-x_{[i;j)}^{-1})} - \sum_{s=i}^{j-1} \frac {\rho_s}{1-x_{[i;j)}^{-1}} \right] \stackrel{\eqref{eqn:mich}}=
$$
\begin{equation}
\label{eqn:has 7}
\stackrel{\eqref{eqn:mich}}= \delta^{-1}_n \left( \sum_{k=1}^\infty \frac {nk}{1-x^{-k}} - \sum^{i < j}_{(i,j) \in \BZ^2/(n,n)\BZ}\left[ \frac {1}{(1-x_{[i;j)})(1-x_{[i;j)}^{-1})} + \sum_{s=i}^{j-1} \frac {\rho_s}{1-x_{[i;j)}^{-1}} \right] \right) \qquad 
\end{equation} 
We leave it as an exercise to the interested reader to use formulas \eqref{eqn:has 1}--\eqref{eqn:has 7} and the Leibniz rule in order to show that \eqref{eqn:leibniz} is equivalent to:
$$
\left[ \CH_m - \sum_{i=1}^n b_i \left( b_i - b_{i-1} + \frac pn \right) \right] Y_m = 0
$$
The only non-trivial part of the computation is the following identity:
\begin{equation}
\label{eqn:sketch}
\sum^{i < j}_{(i,j) \in \BZ^2/(n,n)\BZ} \left[ \frac {j-i}{1 - x_{[i;j)}^{-1}} - \frac 1{(1-x_{[i;j)})(1-x_{[i;j)}^{-1})} - \sum_{s=i+1}^{j} \frac {\rho_s}{1-x_{[i;j)}^{-1}} \right] = 0 \qquad 
\end{equation}
whose proof is just manipulation, and we will now sketch. By definition, we have:
$$
\sum^{i < j}_{(i,j) \in \BZ^2/(n,n)\BZ} \sum_{s=i+1}^{j} \frac {\rho_s}{1-x_{[i;j)}^{-1}} = \mathop{\sum^{a < b}_{(i,j) \in \BZ^2/(n,n)\BZ}}_{(a,b) \in \BZ^2/(n,n)\BZ}^{i < j} \frac {\sum_{s=a}^{b-1} \left(\delta_{\bar{s}}^{\bari} - \delta_{\bar{s}}^{\barj} \right)}{(1-x_{[i;j)}^{-1})(1-x_{[a;b)}^{-1})} = 
$$
\begin{equation}
\label{eqn:sum}
= \frac 12 \mathop{\sum^{a < b}_{(i,j) \in \BZ^2/(n,n)\BZ}}_{(a,b) \in \BZ^2/(n,n)\BZ}^{i < j} \frac {\langle[i;j), [a;b)\rangle + \langle[a;b), [i;j)\rangle}{(1-x_{[i;j)}^{-1})(1-x_{[a;b)}^{-1})}
\end{equation}
where the latter equality holds by symmetrizing $(i,j) \leftrightarrow (a,b)$. It is elementary to see that the numerator of the fraction above is 0 unless one of the numbers $i,j$ is $\equiv$ modulo $n$ to one of the numbers $a,b$. Apart from the terms corresponding to $(i,j) = (a,b)$ modulo $(n,n)\BZ$ (which take care of the middle summand in \eqref{eqn:sketch}), formula \eqref{eqn:sum} is a linear combination of the terms:
$$
\frac 1{(1-x_{[i;j)}^{-1})(1-x_{[i;k)}^{-1})}, \qquad \frac 1{(1-x_{[i;k)}^{-1})(1-x_{[j;k)}^{-1})}, \qquad \frac {-1}{(1-x_{[i;j)}^{-1})(1-x_{[j;k)}^{-1})}
$$
over all $i<j<k$, with simple integer coefficients that the interested reader can easily deduce. It is easy to observe that the sum of the three terms above equals: 
$$
\frac 1{1-x_{[i;k)}^{-1}}
$$
which leads to \eqref{eqn:sketch}. 


\end{proof}

\section{The Proof of geometric relations}
\label{sec:inter}

\subsection{}\label{sub:tautological}

A crucial element in the proof of Theorem \ref{thm:geom} is the ability to calculate the operators \eqref{eqn:fine ops}--\eqref{eqn:smooth ops} in terms of tautological classes, which we will now review (see \cite{Aff} for details). There exist tautological rank $d_i$ vector bundles:
$$
\xymatrix{\CV_i \ar@{.>>}[d] \\ \CM_\bd}
$$
for all $i \in \{1,...,n\}$. Let us consider the ring:
\begin{equation}
\label{eqn:ring}
\Lambda = \BC(u_1,...,u_n,q,\oq)[x_{ia}^{\pm 1}]_{1 \leq i \leq n}^{a \in \BN, \sym}
\end{equation}
where the word $\sym$ refers to functions which are symmetric in $x_{i1},x_{i2},...$ for each $i$ separately. Then we can define a ring homomorphism $\Lambda \longrightarrow K_{\bd}$ for all $\bd \in \nn$:
\begin{equation}
\label{eqn:tautological}
f(...,x_{i1},x_{i2},...) \quad \mapsto \quad \bar{f}_{\bd} :=  f(...,v_{i1},...,v_{id_i},0,......)
\end{equation}
where the symbols $v_{ia}$ are defined by the identity $[\CV_i] = v_{i1}+...+v_{id_i}$ formally in $K_\bd$. We will often write $\bar{f} \in K$ for the image of a symmetric Laurent polynomial under the map \eqref{eqn:tautological}, when the particular degree $\bd \in \nn$ will not be important. \\

\subsection{}\label{sub:variables}

It is often easier to collect variables into a single alphabet:
$$
X = X_1 + ... + X_n = \sum^{a\in \BN}_{1 \leq i \leq n} x_{ia} 
$$
so we may write $f(X)$ instead of $f(...,x_{ia},...)$. Given two alphabets of variables $X$ and $Z$, we may consider the {\bf plethysm} ring homomorphisms:
\begin{equation}
\label{eqn:pleth}
\Lambda \longrightarrow \Lambda[Z, Z^{- 1}]
\end{equation}
$$
f(X) \mapsto f(X \pm Z)
$$
defined on the colored power-sum functions that generate $\Lambda$ as:
$$
f(X) = \sum^{a\in \BN}_{1 \leq i \leq n} x_{ia}^d \quad \mapsto \quad f(X \pm Z) = \sum^{a\in \BN}_{1 \leq i \leq n} x_{ia}^d \pm \sum^{a\in \BN}_{1 \leq i \leq n} z_{ia}^d
$$ 
for all $d \in \BZ \backslash 0$. We will always consider our variables to be {\bf colored} by integers, under the convention $\col x_{ia} = i$ and the almost periodicity property \eqref{eqn:almost}. This property states that we allow ourselves the freedom to replace the variable:
\begin{equation}
\label{eqn:almost}
\Big( x \text{ of color } i-n \Big) \qquad \text{by} \qquad \Big( x\oq^2 \text{ of color } i \Big)
\end{equation}
in all our formulas. The color is multiplicative in the variables, i.e.
$$
\col xx' = \col x + \col x' \qquad \col \frac x{x'} = \col x - \col x'
$$
With this in mind, the rational functions that will appear in this paper may depend on the colors of their variables, the primary example being:
\begin{equation}
\label{eqn:defzeta}
\zeta(x) := \left( \frac {xq - q^{-1}}{x-1} \right)^{\delta_{\col x}^0 - \delta_{\col x + 1}^0}
\end{equation}
whenever $\col x \in [-\frac n2, \frac n2)$. When the colors are arbitrary integers, we may extend the definition of $\zeta$ according to the almost periodicity property \eqref{eqn:almost}, so explicitly:
\begin{equation}
\label{eqn:zeta}
\zeta\left( \frac zw \right) := \left( \frac {zq \oq^{2\left \lceil \frac {i-j}n \right \rceil} - wq^{-1}}{z \oq^{2\left \lceil \frac {i-j}n \right \rceil} - w} \right)^{\delta^0_{i - j \text{ mod }n} - \delta^0_{i - j + 1 \text{ mod }n}}
\end{equation}
if $\col z = i$, $\col w = j$. Note the identity:
\begin{equation}
\label{eqn:identity zeta}
\zeta\left(x \text{ of color } i \right) = \zeta \left(\frac 1{xq^2} \text{ of color }-1-i \right)
\end{equation}
We may extend the notation \eqref{eqn:zeta} multiplicatively to any alphabets of variables $Z = \sum z_{ia}$ and $W = \sum w_{jb}$:
$$
\zeta\left(\frac ZW \right) = \prod_{1 \leq i,j \leq n}^{a,b} \zeta\left( \frac {z_{ia}}{w_{jb}} \right)
$$
If $Z=W$, then we must define instead:
$$
\zeta\left(\frac ZZ \right) = \prod_{(i,a) \neq (j,b)} \zeta\left( \frac {z_{ia}}{z_{jb}} \right)
$$
on account of the fact that $\zeta(x)$ has a pole at $x=1$ if $\col x = 0$. \\

\subsection{}\label{sub:main} 

Tautological classes give a good language for reformulating the operators of Theorem \ref{thm:geom}. For an alphabet $Z = \sum_{1 \leq i \leq n} z_{ia}$, define:
\begin{align}
&\tau_+(Z) = \prod_{1 \leq i \leq n}^{a} \left(\frac {u_{i+1}}q - \frac {z_{ia} q}{u_{i+1}} \right) \label{eqn:deftau 1} \\ &\tau_-(Z) = \prod_{1 \leq i \leq n}^a \left(u_i - \frac {z_{ia}}{u_i} \right) \label{eqn:deftau 2}
\end{align}
Let $Dz = \frac {dz}{2\pi i z}$. The following formulas were proved in \cite{Aff}: \\

\begin{proposition}
\label{prop:1}

For all $i<j \in \BZ$ and any Laurent polynomial $M(z_i,...,z_{j-1})$ with coefficients in $p^{+*}(K_T(\CM_{\bd^+}))$, we have the following equality in $K_T(\CM_{\bd^+})$:
\begin{multline}
p^+_*\Big( M \left(\CL_i,...,\CL_{j-1}\right) [\fZ_{[i;j)}^+] \Big) = \int_{X \prec z_{j-1} \prec ... \prec z_i \prec \{0,\infty\}} \\  \frac {M(z_i,...,z_{j-1}) \prod_{a=i}^{j-1} \left[ \overline{\zeta\left(\frac {z_a}{X} \right)} \tau_+(z_a)  \right] Dz_a}{(q^{-1} - q) \prod_{a=i+1}^{j-1} \left(1 - \frac {z_a}{z_{a-1}q^2}\right)\prod_{i \leq a < b < j} \zeta \left( \frac {z_a}{z_b} \right)} \label{eqn:push 1}
\end{multline}
and if $M$ has coefficients in $p^{-*}(K_T(\CM_{\bd^-}))$, then we have in $K_T(\CM_{\bd^-})$:
\begin{multline}
p^-_*\Big( M \left(\CL_i,...,\CL_{j-1}\right) [\fZ_{[i;j)}^+] \Big) = q^{(d_{i-1}^- - d_{j-1}^-) - (d_i^+ - d_j^+) + i - j}  (-1)^{j-i}  \frac {u_j}{u_i}   \\ \int_{\{0,\infty\} \prec z_{j-1} \prec ... \prec z_i \prec X} \frac {z_i...z_{j-1} M(z_i,...,z_{j-1}) \prod_{a=i}^{j-1} \left[ \overline{\zeta\left(\frac X{z_a} \right)} \tau_-(z_a)  \right]^{-1} Dz_a}{(q^{-1} - q) \prod_{a=i+1}^{j-1} \left(1 - \frac {z_a}{z_{a-1}q^2}\right)\prod_{i \leq a < b < j} \zeta \left( \frac {z_a}{z_b} \right)} \label{eqn:push 2} \qquad
\end{multline}
where $X \prec z_{j-1} ... \prec z_i \prec Y$ means that we integrate the variables $z_{j-1},...,z_i$ over contours contained between the sets $X,Y \subset \BC$, very far away from each other, and ordered with $z_{j-1}$ closest to the set $X$ and $z_i$ closest to the set $Y$. The orientation of the contours is such that the residues in $Y$ are picked up with positive sign. \\

\end{proposition}


\begin{proposition}
\label{prop:2}

For all $i<j \in \BZ$ and any Laurent polynomial $M(z_i,...,z_{j-1})$ with coefficients in $\op^{-*}(K_T(\CM_{\bd^-}))$, we have the following equality in $K_T(\CM_{\bd^-})$:
\begin{multline}
\op^-_*\Big( M \left(\CL_i,...,\CL_{j-1}\right) [\ofZ_{[i;j)}^-] \Big) =  \sum^{t\geq 1}_{i = k_0 < k_1 < ... < k_t = j} \int_{X \prec z_{j-1} = ... = z_{k_{t-1}}  \prec ... \prec z_{k_1-1} = ... = z_{i} \prec \{0,\infty\}} \\ (-1)^{j-i-t}  \frac { M(z_i,...,z_{j-1}) \prod_{a=i}^{j-1} \left[ \overline{\zeta\left(\frac {X}{z_a} \right)} \tau_-(z_a) \right]^{-1} Dz_a }{q^{j-i} (1 - q^{-2}) \prod_{s=1}^{t-1} \left(1 - \frac {z_{k_s-1}}{z_{k_s}} \right)\prod_{i \leq a < b < j} \zeta \left( \frac {z_b}{z_a} \right)} \label{eqn:push 3}
\end{multline}	
and if $M$ has coefficients in $\op^{+*}(K_T(\CM_{\bd^+}))$, then we have in $K_T(\CM_{\bd^+})$:
\begin{multline}
\op^+_*\Big( M \left(\CL_i,...,\CL_{j-1}\right) [\ofZ_{[i;j)}^-] \Big) =  q^{(d_i^+ - d_j^+) - (d_{i-1}^- - d_{j-1}^-) + j-i} (-1)^{j-i} \frac {u_i}{u_j} \\ \sum^{t\geq 1}_{i = k_0 < k_1 < ... < k_t = j} \int_{\{0,\infty\} \prec z_{j-1} = ... = z_{k_{t-1}} \prec ... \prec z_{k_1 - 1} = ... = z_{i} \prec X} \\ (-1)^{j-i-t} \frac {1}{z_i ... z_{j-1}} \frac {M(z_i,...,z_{j-1}) \prod_{a=i}^{j-1} \left[ \overline{\zeta\left(\frac {z_a}{X} \right)} \tau_+(z_a) \right] Dz_a }{q^{j-i}(1 - q^{-2}) \prod_{s=1}^{t-1} \left(1 - \frac {z_{k_s-1}}{z_{k_s}} \right) \prod_{i \leq a < b < j} \zeta \left( \frac {z_b}{z_a} \right)}  \label{eqn:push 4}
\end{multline}
The right-hand side of \eqref{eqn:push 3} should be interpreted as follows: every summand corresponds to a way to divide the variables $z_{j-1},...,z_i$ into $t$ consecutive groups. Set all the variables in the $i$--th group equal to some variable $y_a$, and then integrate the resulting function along the contours $X \prec y_t \prec ... \prec y_1 \prec \{0,\infty\}$ (the notation $\prec$ is explained in Proposition \ref{prop:1}). The integral in \eqref{eqn:push 4} is defined analogously. \\

\end{proposition}

\begin{proposition}
\label{prop:smooth}

For all $k \in \BN$ and any symmetric Laurent polynomial $M$ in the tautological vector bundles $\CL_1,...,\CL_n$, we have the following equality in $K_T(\CM_{\bd^+})$:
\begin{multline}
\pi^+_*\Big( M(\CL_1,...,\CL_n) \cdot [\fW^+_k] \Big) = \sum^{\text{compositions}}_{k_1+...+k_t = k} \frac 1{k_1(k_1+k_2)...(k_1+...+k_t)} \\ \int_{X \prec y_1 \prec ... \prec y_t \prec \{0,\infty\}} M(Y) \cdot \eta\left( \frac YY \right) \left[ \overline{\zeta \left( \frac YX \right)} \tau_+(Y) \right] \prod_{s=1}^t Dy_s \label{eqn:push 5}
\end{multline}
and the following equality in $K_T(\CM_{\bd^-})$:
\begin{multline}
\pi^-_*\Big( M(\CL_1,...,\CL_n) \cdot [\fW^-_k] \Big) = \sum^{\text{compositions}}_{k_1+...+k_t = k} \frac {1}{k_1(k_1+k_2)...(k_1+...+k_t)} \\ \int_{\{0, \infty\} \prec y_1 \prec ... \prec y_t \prec X} M(Y) \cdot \eta\left( \frac YY \right) \left[ \overline{\zeta \left( \frac XY \right)} \tau_-(Y) \right]^{-1} \prod_{s=1}^t Dy_s \label{eqn:push 6}
\end{multline}
where in the formulas above, $Y = Y_1+...+Y_n$ with $Y_i = \sum_{s=1}^t \sum_{a=0}^{k_s-1} y_s \oq^{-2a}$, and:
$$
\eta \left( \frac {\sum_i z_i}{\sum_j w_j} \right) = \prod_{i,j} \eta \left( \frac {z_i}{w_j} \right), \qquad \text{where} \qquad \eta(x) = (1-x)^{\delta_{\col x}^0 - \delta_{\col x}^{-1}}
$$
The integrals in \eqref{eqn:push 5}--\eqref{eqn:push 6} are interpreted as in Propositions \ref{prop:1} and \ref{prop:2}. \\

\end{proposition}

\begin{remark}
\label{rem:poles}

Although not explicitly spelled out in \cite{Aff}, there are immediate analogues of formulas \eqref{eqn:push 1}--\eqref{eqn:push 3} when the Laurent polynomial $M(z_i,...,z_{j-1})$ is replaced by a rational function with poles in some set $X' \sqcup \{0,\infty\}$. The only modification the reader needs to make in order for formulas \eqref{eqn:push 1}--\eqref{eqn:push 6} is to change $\{0,\infty\}$ by $X' \sqcup \{0,\infty\}$ in the subscripts of all the integrals. We leave the proof of this more general, but straightforward, fact to the interested reader. \\  

\end{remark}

\subsection{}\label{sub:e tautological}


Tautological classes are very useful since they allow us to compute the vector bundle $E$ of Subsection \ref{sub:e}. For example, the following result was proved in \cite{Aff}:
\begin{equation}
\label{eqn:ke}
[\CE] = \sum_{i=1}^{n} \left[\left(1-\frac 1{q^2}\right)\left(\frac {\CV_i}{\CV'_{i-1}} - \frac {\CV_i}{\CV'_i} \right) + \frac {u_{i+1}^2}{\CV'_iq^2} + \frac {\CV_i}{u_i^2}\right] \in K_\bd \otimes K_{\bd'}
\end{equation}
where $\CV_i$ and $\CV_i'$ are the pull-backs of tautological bundles from the factors $\CM_\bd$ and $\CM_{\bd'}$, respectively. To keep our formulas simple, in \eqref{eqn:ke} and beyond, we will write $\CV$ instead of $[\CV]$ when computing $K-$theory classes such as \eqref{eqn:ke}, as well as:
$$
\frac {\CV'}{\CV} \qquad \text{instead of} \qquad [\CV'] \cdot [\CV^\vee]
$$
As a consequence of \eqref{eqn:ke}, we conclude that:
$$
[\wedge^\bullet(\CE^\vee,m) ] = \prod_{i=1}^{n} \frac {\left(1 - \frac {\CV_{i-1}' q^{2m}}{\CV_i} \right)\left(1 - \frac {\CV_{i}'q^{2m+2}}{\CV_i} \right)}{\left(1 - \frac {\CV_{i-1}' q^{2m+2}}{\CV_i} \right)\left(1 - \frac {\CV'_{i}q^{2m}}{\CV_i} \right)} \left(1 - \frac {\CV_i' q^{2m+2}}{u_{i+1}^2}\right) \left(1 - \frac {u_i^2q^{2m}}{\CV_i}\right) 
$$
Using the renormalization in \eqref{eqn:modification}, we conclude that:
\begin{equation}
\label{eqn:panayiotou} 
[\widetilde{\wedge}^\bullet (\CE^\vee,m)] = \overline {\zeta\left(\frac {X' q^{2m}}{X} \right) \tau_+ \left( X' q^{2m} \right) \tau_-\left( X q^{-2m} \right)} \cdot q^{m(|\bd| - |\bd'|)} \lambda_{\bd'}
\end{equation}
where $X$, $X'$ are place-holders for the tautological vector bundles $\{\CV_i\}$, $\{\CV_i'\}$ on $\CM_\bd$, $\CM_{\bd'}$, respectively. Formula \eqref{eqn:panayiotou} gives a class $\in K_\bd \otimes K_{\bd'}$ via \eqref{eqn:tautological}. \\


\begin{proof} {\bf of Proposition \ref{prop:comm root e}:} Consider the following diagrams:
\begin{equation}
\label{eqn:big diagram 1}
\xymatrix{ & & \CM_{\bd} \times \CM_{\bd'} \ar@/_3pc/[llddd]_{\pi_1} \ar@/^3pc/[rrddd]^{\pi_2} & & \\
& & \CM_{\bd} \times \fZ_{[i;j)} \ar[ld]_{\id \times p_+} \ar[rd]   \ar[u]^{\id \times p_-} & & \\
& \CM_{\bd} \times \CM_{\bd'+[i;j)} \ar[ld] \ar[rd] & & \fZ_{[i;j)} \ar[ld]^{p_+} \ar[rd]_{p_-} & \\ 
\CM_\bd & & \CM_{\bd'+[i;j)} & & \CM_{\bd'}}
\end{equation}
\begin{equation}
\label{eqn:big diagram 2}
\xymatrix{ & & \CM_{\bd} \times \CM_{\bd'}  \ar@/_3pc/[llddd]_{\pi_1'} \ar@/^3pc/[rrddd]^{\pi'_2} & & \\
& & \fZ_{[i;j)} \times \CM_{\bd'} \ar[ld] \ar[rd]^{p_-' \times \id}   \ar[u]_{p_+' \times \id} & & \\
& \fZ_{[i;j)} \ar[ld]^{p'_+} \ar[rd]_{p'_-} & & \CM_{\bd-[i;j)} \times \CM_{\bd'} \ar[ld] \ar[rd] & \\ 
\CM_{\bd} & & \CM_{\bd-[i;j)} & & \CM_{\bd'}}
\end{equation}	
Because of \eqref{eqn:fine ops}, we have:
\begin{align}
A_m e_{[i;j)} &= \pi_{1*} \Big(\Upsilon_{[i;j)} \cdot \pi_2^* \Big) \label{eqn:first} \\
e_{[i;j)} A_m &= \pi'_{1*} \Big(\Upsilon'_{[i;j)} \cdot {\pi_2'}^{*} \Big) \label{eqn:second} 
\end{align}
where:
\begin{align*}
&\Upsilon_{[i;j)} = (\id \times p_-)_* \Big[ [\fZ_{[i;j)}^+] \cdot \widetilde{\wedge}^\bullet ((\id \times p_+)^*(\CE^\vee),m) \Big]  \\
&\Upsilon'_{[i;j)} = (p_+' \times \id)_* \Big[ [\fZ_{[i;j)}^+] \cdot \widetilde{\wedge}^\bullet ((p'_- \times \id)^*(\CE^\vee),m) \Big]
\end{align*}
As a consequence of \eqref{eqn:panayiotou}, it is easy to check the following formulas in $K$--theory:
\begin{align*}
&\widetilde{\wedge}^\bullet ((\id \times p_+)^*\CE^\vee,m) = \widetilde{\wedge}^\bullet ((\id \times p_-)^*\CE^\vee,m) \prod_{a=i}^{j-1} \left[ \zeta \left(\frac {\CL_a q^m}{X} \right) \tau_+\left( \CL_a q^m \right) \right] \cdot q^{m(i-j)} \frac {\lambda_{\bd'+[i;j)}}{\lambda_{\bd'}} \\ 
& \widetilde{\wedge}^\bullet ((p'_- \times \id)^*\CE^\vee,m) =  \widetilde{\wedge}^\bullet ((p'_+ \times \id)^*\CE^\vee,m) \prod_{a=i}^{j-1} \left[ \zeta \left( \frac {X' q^m}{\CL_a} \right) \tau_-(\CL_a q^{-m}) \right]^{-1} \cdot q^{m(j-i)}
\end{align*}
where $\frac {\lambda_{\bd'+[i;j)}}{\lambda_{\bd'}} = (-1)^{j-i} \frac {u_i}{u_j} \frac {q^{j-i + \langle \bd', [ i;j) \rangle + \langle [i;j), \bd' + [i;j) \rangle}}{\CL_i... \CL_{j-1}}$. Therefore, we conclude that:	
\begin{align*}
&\Upsilon_{[i;j)} = \widetilde{\wedge}^\bullet (\CE^\vee,m) \cdot (\id \times p_-)_* \left( [\fZ_{[i;j)}^+] \cdot \prod_{a=i}^{j-1} \left[ \zeta \left(\frac {\CL_a q^m}{X} \right) \tau_+\left( \CL_a q^m \right) \right] \right) \cdot q^{m(i-j)} \cdot \frac {\lambda_{\bd'+[i;j)}}{\lambda_{\bd'}} \\
&\Upsilon'_{[i;j)} = \widetilde{\wedge}^\bullet (\CE^\vee,m) \cdot (p_+' \times \id)_* \left( [\fZ_{[i;j)}^+] \cdot  \prod_{a=i}^{j-1} \left[ \zeta \left( \frac {X' q^m}{\CL_a} \right) \tau_-(\CL_a q^{-m}) \right]^{-1} \right) \cdot q^{m(i-j)}
\end{align*}
Using formulas \eqref{eqn:push 1} and \eqref{eqn:push 2} (with Remark \ref{rem:poles}), we obtain:
\begin{align}
&\Upsilon_{[i;j)} = \widetilde{\wedge}^\bullet (\CE^\vee,m) \int_{X \sqcup \{0,\infty\} \prec z_{j-1} \prec ... \prec z_i \prec X'} \frac {Q_{[i;j)}(z_iq^{2m},...,z_{j-1} q^{2m})}{q^{-1} - q} \prod_{a=i}^{j-1} Dz_a \label{eqn:upsilon 1} \\
&\Upsilon_{[i;j)}' = \widetilde{\wedge}^\bullet (\CE^\vee,m) \int_{X \prec z_{j-1} \prec ... \prec z_i \prec X' \sqcup \{0,\infty\}} \frac {Q_{[i;j)}(z_i,...,z_{j-1})}{q^{-1} - q} \prod_{a=i}^{j-1} Dz_a \label{eqn:upsilon 2}
\end{align}
where:
$$
Q_{[i;j)}(z_i,...,z_{j-1}) = \prod_{a=i}^{j-1} \left[ \overline{\frac {\zeta\left(\frac {z_a}{X} \right)}{\zeta \left( \frac {X' q^m}{z_a} \right)}} \frac {\tau_+(z_a)}{\tau_-(z_aq^{-m})}  \right] \frac {q^{m(i-j)}}{\prod_{a=i+1}^{j-1} \left(1 - \frac {z_a}{z_{a-1}q^2}\right)\prod_{i \leq a < b < j} \zeta \left( \frac {z_a}{z_b} \right)}
$$ 
Observe that the only difference between the formulas \eqref{eqn:upsilon 1} and \eqref{eqn:upsilon 2} is that the contours of integration are different. Therefore, the difference between \eqref{eqn:upsilon 1} and \eqref{eqn:upsilon 2} stems from the residues at 0 and at $\infty$ of the function $Q_{[i;j)}$. It is easy to see that the only non-zero such residues are:
\begin{multline*}
\underset{z_{i} = \infty}{\text{Res}} \frac {Q_{[i;j)}(z_i,...,z_{j-1})}{z_i} = \\ = Q_{[i+1;j)}(z_{i+1},...,z_{j-1}) \cdot q^{d_i - d_{i+1} + d_i' - d_{i-1}'} \frac {q^{1+m}u_i}{u_{i+1}} q^{- \langle [i;i+1), [i+1;j) \rangle}
\end{multline*}
\begin{multline*} 
\underset{z_{j-1} = 0}{\text{Res}} \frac {Q_{[i;j)}(z_i,...,z_{j-1})}{z_i} = \\ = Q_{[i;j-1)}(z_{i},...,z_{j-2}) \cdot q^{-d_{j-1} + d_{j} + d_{j-2}' - d_{j-1}'} \frac {q^{-m-1}u_j}{u_{j-1}} q^{- \langle [i;j-1), [j-1;j) \rangle}
\end{multline*}
Then the residue theorem reads:
\begin{multline*}
\Upsilon'_{[i;j)} + \Upsilon'_{[i+1;j)} q^{d_i - d_{i+1} + d_i' - d_{i-1}'} \frac {q^{1+m}u_i}{u_{i+1}} q^{\langle [i;j-1), [j-1;j) \rangle} = \\ = \Upsilon_{[i;j)} + \Upsilon_{[i;j-1)} q^{-d_{j-1} + d_{j} + d_{j-2}' - d_{j-1}'} \frac {q^{-m-1}u_j}{u_{j-1}} q^{\langle [i;i+1), [i+1;j) \rangle}
\end{multline*}
Using \eqref{eqn:first}, \eqref{eqn:second} and \eqref{eqn:psi}, the equality above is equivalent to \eqref{eqn:comm1}. \\

\end{proof}

\begin{proof} \textbf{of Proposition \ref{prop:comm root f}:} Consider the following diagrams:
\begin{equation}
\label{eqn:big diagram 3}
\xymatrix{ & & \CM_{\bd} \times \CM_{\bd'} \ar@/_3pc/[llddd]_{\pi_1} \ar@/^3pc/[rrddd]^{\pi_2} & & \\
& & \CM_{\bd} \times \ofZ_{[i;j)} \ar[ld]_{\id \times \op_-} \ar[rd]   \ar[u]^{\id \times \op_+} & & \\
& \CM_{\bd} \times \CM_{\bd'-[i;j)} \ar[ld] \ar[rd] & & \ofZ_{[i;j)} \ar[ld]^{\op_-} \ar[rd]_{\op_+} & \\ 
\CM_\bd & & \CM_{\bd'-[i;j)} & & \CM_{\bd'}}
\end{equation}
\begin{equation}
\label{eqn:big diagram 4}
\xymatrix{ & & \CM_{\bd} \times \CM_{\bd'}  \ar@/_3pc/[llddd]_{\pi_1'} \ar@/^3pc/[rrddd]^{\pi'_2} & & \\
& & \ofZ_{[i;j)} \times \CM_{\bd'} \ar[ld] \ar[rd]^{\op_+' \times \id}   \ar[u]_{\op_-' \times \id} & & \\
& \ofZ_{[i;j)} \ar[ld]^{\op'_-} \ar[rd]_{\op'_+} & & \CM_{\bd+[i;j)} \times \CM_{\bd'} \ar[ld] \ar[rd] & \\ 
\CM_{\bd} & & \CM_{\bd+[i;j)} & & \CM_{\bd'}}
\end{equation}	
Because of \eqref{eqn:eccentric ops}, we have:
\begin{align}
A_m f_{[i;j)} &= \pi_{1*} \Big(\Upsilon_{[i;j)} \cdot \pi_2^* \Big) \label{eqn:third} \\
f_{[i;j)} A_m &= \pi'_{1*} \Big(\Upsilon'_{[i;j)} \cdot {\pi_2'}^{*} \Big) \label{eqn:fourth}
\end{align}
where:
\begin{align*}
&\Upsilon_{[i;j)} = (\id \times \op_+)_* \Big[ [\ofZ_{[i;j)}^-] \cdot \widetilde{\wedge}^\bullet ((\id \times \op_-)^*(\CE^\vee),m) \Big]  \\
&\Upsilon'_{[i;j)} = (\op_-'  \times \id)_* \Big[ [\ofZ_{[i;j)}^-] \cdot \widetilde{\wedge}^\bullet ((\op'_+ \times \id)^*(\CE^\vee),m) \Big]
\end{align*}
As a consequence of \eqref{eqn:panayiotou}, it is easy to check the following formulas in $K$--theory:
\begin{align*}
&\widetilde{\wedge}^\bullet ((\id \times \op_-)^*\CE^\vee,m) = \widetilde{\wedge}^\bullet ((\id \times \op_+)^*\CE^\vee,m) \prod_{a=i}^{j-1} \left[ \zeta \left( \frac {\CL_a q^{2m}}{X} \right) \tau_+ \left( \CL_a q^{2m} \right) \right]^{-1} \cdot q^{m(j-i)} \frac {\lambda_{\bd' - [i;j)}}{\lambda_{\bd'}} \\
&\widetilde{\wedge}^\bullet ((\op'_+ \times \id)^*\CE^\vee,m) = \widetilde{\wedge}^\bullet ((\op'_- \times \id)^*\CE^\vee,m) \prod_{a=i}^{j-1} \left[ \zeta \left( \frac {X' q^{2m}}{\CL_a} \right) \tau_- \left(\CL_a q^{-2m} \right) \right] \cdot q^{m(j-i)}
\end{align*}
where $\frac {\lambda_{\bd'-[i;j)}}{\lambda_{\bd'}} = (-1)^{j-i} \frac{u_j}{u_i} q^{i-j - \langle \bd' - [i;j), [ i;j) \rangle - \langle [i;j), \bd' \rangle} \CL_i... \CL_{j-1}$. Therefore:
\begin{align*}
&\Upsilon_{[i;j)} = \widetilde{\wedge}^\bullet (\CE^\vee,m) \cdot (\id \times \op_+)_* \left( [\ofZ_{[i;j)}^-] \cdot \prod_{a=i}^{j-1} \left[ \zeta \left( \frac {\CL_a q^{2m}}{X} \right) \tau_+ \left( \CL_a q^{2m} \right) \right]^{-1} \right)  \cdot q^{m(j-i)} \frac {\lambda_{\bd' - [i;j)}}{\lambda_{\bd'}} \\
&\Upsilon'_{[i;j)} = \widetilde{\wedge}^\bullet (\CE^\vee,m) \cdot (\op_-'  \times \id)_* \left( [\ofZ_{[i;j)}^-] \cdot \prod_{a=i}^{j-1} \left[ \zeta \left( \frac {X' q^{2m}}{\CL_a} \right) \tau_- \left(\CL_a q^{-2m} \right) \right] \right) \cdot q^{m(j-i)}
\end{align*}
Using formulas \eqref{eqn:push 3} and \eqref{eqn:push 4} (with Remark \ref{rem:poles}), we obtain:
\begin{align} 
&\Upsilon_{[i;j)} = \widetilde{\wedge}^\bullet (\CE^\vee,m) \sum^{t\geq 1}_{i = k_0 < k_1 < ... < k_t = j} (-1)^{j-i-t} V_{k_0,...,k_t} \label{eqn:upsilon 3} \\
&\Upsilon'_{[i;j)} = \widetilde{\wedge}^\bullet (\CE^\vee,m) \sum^{t\geq 1}_{i = k_0 < k_1 < ... < k_t = j} (-1)^{j-i-t} V_{k_0,...,k_t}' \label{eqn:upsilon 4}
\end{align}
where:
\begin{align*}
&V_{k_0,...,k_t}  = \int_{X \sqcup \{0,\infty\} \prec z_{j-1} = ... = z_{k_{t-1}}  \prec ... \prec z_{k_1-1} = ... = z_{i} \prec  X'} \frac {Q_{[i;j)}(z_i,...,z_{j-1})}{1 - q^{-2}} \prod_{a=i}^{j-1} Dz_a \\
&V_{k_0,...,k_t}'  = \int_{X \prec z_{j-1} = ... = z_{k_{t-1}}  \prec ... \prec z_{k_1-1} = ... = z_{i} \prec X' \sqcup \{0,\infty\}} \frac {Q'_{[i;j)}(z_i,...,z_{j-1})}{1 - q^{-2}} \prod_{a=i}^{j-1} Dz_a
\end{align*}
and:
\begin{align*}
&Q_{[i;j)}(z_i,...,z_{j-1}) = \prod_{a=i}^{j-1} \left[ \overline{\frac {\zeta \left( \frac {z_a}{X'} \right)}{\zeta \left( \frac {z_a q^{2m} }{X} \right)}} \frac {\tau_+(z_a)}{\tau_+(z_a q^{2m})} \right] \frac {q^{(m-1)(j-i)}}{\prod_{s=1}^{t-1} \left(1 - \frac {z_{k_s-1}}{z_{k_s}} \right) \prod_{i \leq a < b < j} \zeta \left( \frac {z_b}{z_a} \right)} \\
&Q_{[i;j)}'(z_i,...,z_{j-1}) = \prod_{a=i}^{j-1} \left[ \overline{\frac {\zeta \left( \frac {X' q^{2m}}{z_a} \right)}{\zeta \left( \frac {X}{z_a} \right)}} \frac {\tau_-(z_a q^{-2m})}{\tau_-(z_a)} \right] \frac {q^{(m-1)(j-i)}}{\prod_{s=1}^{t-1} \left(1 - \frac {z_{k_s-1}}{z_{k_s}} \right) \prod_{i \leq a < b < j} \zeta \left( \frac {z_b}{z_a} \right)}
\end{align*}
Because of the identity \eqref{eqn:identity zeta} , one observes that:
$$
Q_{[i;j)}(z_i,...,z_{j-1}) = Q'_{[i+1;j+1)}(z_iq^{2+2m},...,z_{j-1}q^{2+2m})
$$
Therefore, we conclude that the quantities $V_{k_0,...,k_t}$ and $V_{k_0+1,...,k_t+1}'$ are given by integrals of one and the same function, but the contours differ in the placement of $\{0,\infty\}$. It is easy to see that the only non-zero residues of the function $Q_{[i;j)}$ at either 0 or $\infty$ are:
\begin{multline*}
\underset{z_{j-1} = ... = z_{k_{t-1}} = \infty}{\text{Res}} \frac {Q_{[i;j)}(z_i,...,z_{j-1})}{z_{j-1} ... z_{k_{t-1}}} =  Q_{[i;k_{t-1})}(z_i,...,z_{k_{t-1} - 1}) r_{\left \lfloor \frac {j-k_{t-1}}n \right \rfloor} \\ q^{-(m+1)(j-k_{t-1})+ (d_{k_{t-1}}' - d_j') - (d_{k_{t-1}} - d_j)} q^{- \langle [k_{t-1},j), [i;k_{t-1}) \rangle}
\end{multline*}
\begin{multline*}
\underset{z_{k_1-1} = ... = z_i = 0}{\text{Res}} \frac {Q_{[i;j)}(z_i,...,z_{j-1})}{z_{k_1-1} ... z_i} =  Q_{[k_1;j)}(z_{k_1},...,z_{j-1}) r_{\left \lfloor \frac {k_1-i}n \right \rfloor} \\  q^{(m-1)(k_1-i) + (d_{k_1}' - d_i') - (d_{k_1} - d_i)} q^{- \langle [k_1;j), [i;k_1) \rangle}
\end{multline*}
where $r_{\left \lfloor \frac {j-i}n \right \rfloor} = \prod_{i \leq a < b < j} \zeta\left(\frac {z_b}{z_a} \right)^{-1} \Big|_{z_a \text{ of color } a \mapsto 1} = \prod_{l=1}^{\left \lfloor \frac {j-i}n \right \rfloor} \frac {q \oq^{2l} - q^{-1}}{\oq^{2l}-1}$. Thus:
\begin{multline*}
V'_{k_0+1,...,k_t+1} - V'_{k_1+1,...,k_t+1} r_{\left \lfloor \frac {k_1-i}n \right \rfloor} q^{(m-1)(k_1-i) + (d_{k_1}' - d_i') - (d_{k_1} - d_i)} q^{- \langle [k_1+1;j+1), [i+1;k_1+1) \rangle} \\ = V_{k_0,...,k_t} - V_{k_0,...,k_{t-1}} r_{\left \lfloor \frac {j-k_{t-1}}n \right \rfloor} q^{-(m+1)(j-k_{t-1}) + (d_{k_{t-1}}' - d_j') - (d_{k_{t-1}} - d_j)} q^{- \langle [k_{t-1},j), [i;k_{t-1}) \rangle}
\end{multline*}
Summing up over all $i = k_0 < ... < k_t = j$, we conclude that:
\begin{multline*}
\Upsilon'_{[i+1;j+1)} + \sum_{a=i+1}^j  \Upsilon'_{[a+1;j+1)} (-q)^{i-a} r_{\left \lfloor \frac {a-i}n \right \rfloor}  q^{m(a-i) + (d_a' - d_i') - (d_a - d_i)} q^{- \langle [a+1;j+1), [i+1;a+1) \rangle} \\ = \Upsilon_{[i;j)} + \sum_{a=i}^{j-1} \Upsilon_{[i;a)} (-q)^{a-j} r_{\left \lfloor \frac {j-a}n \right \rfloor}  q^{-m(j-a)+ (d_a' - d_j') - (d_a - d_j)}q^{- \langle [a,j), [i;a) \rangle}
\end{multline*}
Using \eqref{eqn:third}, \eqref{eqn:fourth} and \eqref{eqn:psi}, the equality above is equivalent to \eqref{eqn:comm2}. \\

\end{proof}

\begin{proof} {\bf of Proposition \ref{prop:comm group}:} We will prove \eqref{eqn:comm3}, and leave the analogous formula \eqref{eqn:comm4} to the interested reader. Running the argument in the proof of Propositions \ref{prop:comm root e} and \ref{prop:comm root f} with the modification $\fZ_{[i;j)} \leadsto \fW_k$ implies that: 
\begin{align*}
A_m g_k &= \pi_{1*} \Big(\Upsilon_{[i;j)} \cdot \pi_2^* \Big) \\
g_k A_m &= \pi'_{1*} \Big(\Upsilon'_{[i;j)} \cdot {\pi_2'}^{*} \Big) 
\end{align*}
where:
\begin{align}
&\Upsilon_{k} = \widetilde{\wedge}^\bullet (\CE^\vee,m) \sum^{\text{compositions}}_{k_1+...+k_t = k} \frac {B_{|k_1,...,k_t}}{\#_{k_1,...,k_t}} \label{eqn:upsilon} \\
&\Upsilon_{k}' = \widetilde{\wedge}^\bullet (\CE^\vee,m) \sum^{\text{compositions}}_{k_1+...+k_t = k} \frac {B_{k_1,...,k_t|}}{\#_{k_1,...,k_t}} \label{eqn:upsilon prim}
\end{align}
where $\#_{k_1,...,k_t} = k_1(k_1+k_2)...(k_1+...+k_t)$ and:
$$
B_{k_1,...,k_s|k_{s+1},...k_t} = \oq^{k^2} \int_{X \prec y_1 \prec ... \prec y_s \prec \{0,\infty\} \prec y_{s+1} \prec ... \prec y_t \prec X'} Q_{k_1,...,k_t}(y_1,...,y_t) \prod_{s=1}^t Dy_s \label{eqn:b}
$$
Above, we consider the alphabet $Y = \sum_{i=1}^n \sum_{s=1}^t \sum_{a=0}^{k_s-1} \underbrace{y_s \oq^{-2a}}_{\text{color }i}$, and let:
$$
Q_{k_1,...,k_t}(y_1,...,y_t) = \eta \left( \frac YY \right) \overline{\frac {\zeta\left(\frac {Y}{X} \right)}{\zeta \left( \frac {X' q^{2m}}{Y} \right)}} \frac {\tau_+(Y)}{\tau_-(Yq^{-2m})} \cdot q^{-mnk}
$$
It is easy to see that the residues of the rational function above at $y_s \in \{0,\infty\}$ are:
\begin{align*}
&\underset{y_s = \infty}{\text{Res}} \frac {Q_{k_1,...,k_t}(y_1,...,y_t)}{y_s} = (-1)^{k_s-1} \frac {\oq^{-2k k_s +k_s^2 + k_s}}{1-\oq^{2k_s}} Q_{k_1,...,\widehat{k_s},...,k_t}(y_1,...,\widehat{y_s},...,y_t) q^{(1+m)nk_s} \oq^{k_s} \\ 
&\underset{y_s = 0}{\text{Res}} \frac {Q_{k_1,...,k_t}(y_1,...,y_t)}{y_s} = (-1)^{k_s-1} \frac {\oq^{-2k k_s +k_s^2 + k_s}}{1-\oq^{2k_s}} Q_{k_1,...,\widehat{k_s},...,k_t}(y_1,...,\widehat{y_s},...,y_t) q^{-(m+1)nk_s} \oq^{-k_s}
\end{align*}
which implies that:
$$
B_{...,k_{s-1},k_s|k_{s+1},...} - B_{...,k_{s-1}|k_s,k_{s+1},...} = \frac {(-1)^{k_s-1}}{\oq^{k_s}-\oq^{-k_s}} \left( q^{(\gamma + mn)k_s} - q^{-(mn+\gamma) k_s} \right) B_{...,k_{s-1}|k_{s+1},...}
$$
If we iterate the formula above, we obtain:
$$
B_{k_1,...,k_t|} = \sum^{j \geq 0}_{1 \leq a_1 < ... < a_j \leq t} \prod_{i=1}^j (-1)^{k_{a_i}-1} \frac {q^{(mn + \gamma)k_{a_i}} - q^{-(mn+\gamma) k_{a_i}}}{\oq^{k_{a_i}}-\oq^{-k_{a_i}}} B_{|...,\widehat{k_{a_1}},...,\widehat{k_{a_j}},...}
$$
It is a straightforward exercise in manipulating generating series (which is left to the interested reader), that the formula above implies the identity:
\begin{multline*}
\left(\sum_{k_1,...,k_t \in \BN}^{t \geq 0} \frac {B_{k_1,...,k_t|}}{\#_{k_1,...,k_t}} \cdot x^{k_1+...+k_t} \right) = \left(\sum_{k_1,...,k_t \in \BN}^{t \geq 0} \frac {B_{|k_1,...,k_t}}{\#_{k_1,...,k_t}} \cdot x^{k_1+...+k_t} \right) \cdot \\ \cdot \exp \left(\sum_{k=1}^\infty \frac {(-1)^{k-1} x^k}k \cdot \frac {q^{(\gamma + mn)k} - q^{-(mn+\gamma) k}}{\oq^{k}-\oq^{-k}} \right)
\end{multline*}
Then \eqref{eqn:upsilon} and \eqref{eqn:upsilon prim} imply that:
$$
\left( \sum_{k=0}^\infty g_k x^k \right) A_m = A_m \left( \sum_{k=0}^\infty g_k x^k\right) \exp \left(\sum_{k=1}^\infty \frac {(-1)^{k-1} x^k}k \cdot \frac {q^{(\gamma + mn)k} - q^{-(mn+\gamma) k}}{\oq^{k}-\oq^{-k}} \right)
$$
and by taking the logarithm, we conclude \eqref{eqn:comm3}.

\end{proof}

\end{document}